\documentclass[11pt, leqno]{article}
\usepackage{amssymb,latexsym}

\newcommand{\T}{{\mathbb T}}
\newcommand{\C}{{\mathbb C}}
\newcommand{\R}{{\mathbb R}}

\newcommand{\Q}{{\mathbb Q}}
\newcommand{\Z}{{\mathbb Z}}
\newcommand{\CP}{{\mathbb CP}}

\begin{document}

\title{Symplectic capacities of toric manifolds and related results}

\author{\\Guangcun Lu
\thanks{Partially supported by the NNSF 19971045 and 10371007 of China.}\\
Department of Mathematics,  Beijing Normal University\\
Beijing 100875,   P.R.China\\
(gclu@bnu.edu.cn)}
\date{Received: December 15, 2004\\
Accepted: August 22, 2005} \maketitle \vspace{-0.1in}

\abstract{In this paper we  give concrete estimations for the
pseudo symplectic capacities of toric  manifolds in combinatorial
data.  Some examples are given to show that our estimates can
compute their pseudo symplectic capacities. As applications we
also estimate the symplectic capacities of the polygon spaces.
Other related results are impacts of symplectic blow-up on
symplectic capacities, symplectic packings in symplectic toric
manifolds, the Seshadri constant of an ample line bundle on toric
manifolds, and symplectic capacities of symplectic manifolds with
$S^1$-action.
 \vspace{-0.1in}

\section{Introduction and main results}

The symplectic capacities are the important tools of study  of
symplectic topology. There are several symplectic capacities. The
typical two of them are the Gromov symplectic width ${\cal W}_G$
and Hofer-Zehnder capacity $c_{HZ}$ (cf. [Gr] and [HZ]). However,
for a general symplectic manifold $(M,\omega)$ it is very
difficult to compute or estimate ${\cal W}_G(M,\omega)$ and
$c_{HZ}(M,\omega)$; see [Gin], [Lu3] and reference therein for the
related results. It is well-known that the toric manifolds are a
very beautiful family of K\"ahler manifolds admitting a
combinatorial description. They also are rational and thus
uniruled. So their pseudo symplectic capacities all are finite
(cf. [Lu3]). The main aim of this paper is to estimate their
(pseudo) symplectic capacities in terms of combinatorial data.
Part results were announced in [Lu2] though they should be
restricted to toric Fano manifolds as showed below.\vspace{2mm}

Firstly, we briefly recall the typical {\it pseudo symplectic
capacity} introduced in [Lu3]. For its properties and applications
the reader refer to [Lu3]. Given a connected symplectic manifold
$(M,\omega)$ of dimension $2n$ and a smooth function $H$ on it let
$X_H$ denote the symplectic gradient of $H$. An isolated critical
point $p$ of $H$ is called {\it admissible} if the spectrum of the
linear transformation $DX_H(p): T_pM\to T_pM$ is contained in
$\C\setminus\{\lambda i\,|\, 2\pi\le\pm\lambda<+\infty\}$. For two
given nonzero homology classes $\alpha_0,\alpha_\infty\in
H_\ast(M)$ we denote by
$${\cal H}_{ad}(M,\omega;\alpha_0,\alpha_\infty)\quad{\rm (resp.}\quad
{\cal H}_{ad}^\circ(M,\omega;\alpha_0,\alpha_\infty)\;)$$ the set
of all smooth functions on $M$ for which there exist two smooth
compact submanifolds $P$ and $Q$ of $M$ with connected smooth
boundaries and of codimension zero such that  the following
condition groups (a)(b)(c)(d)(e)(f) (resp.
(a)(b)(c)(d)(e)($f^\prime$)) are satisfied:
\begin{description}
\item[(a)] $P\subset {\it Int}(Q)$ and $Q\subset {\it Int}(M)$;
\item[(b)] $H|_P=0$ and $H|_{M-{\it Int}(Q)}=\max H$;
\item[(c)] $0\le H\le \max H$;
\item[(d)] There exist  chain representatives of $\alpha_0$ and $\alpha_\infty$,
           still denoted by $\alpha_0,\alpha_\infty$, such that
           ${\it supp}(\alpha_0)\subset{\it Int}(P)$ and
           ${\it supp}(\alpha_\infty)\subset M\setminus Q$;
\item[(e)]   There are no critical values in $(0,\varepsilon)\cup
(\max H-\varepsilon, \max H)$    for a small
$\varepsilon=\varepsilon(H)>0$;
\item[(f)] The Hamiltonian system
$\dot x=X_H(x)$ on $M$ has no nonconstant
            periodic solutions of period less than $1$;
\item[(${\bf f}^\prime$)] The Hamiltonian system $\dot x=X_H(x)$ on $M$ has no nonconstant
            contractible periodic solutions of period less than $1$.
\end{description}
If $\alpha_0\in H_0(M)$ can be represented by a point we allow $P$ to be an empty set.
If $M$ is a closed manifold and $\alpha_\infty\in H_0(M)$ is represented by a point,
we also allow $Q=M$.

The pseudo symplectic capacities of Hofer-Zehnder type are defined by
$$\left\{\begin{array}{ll}
C_{HZ}^{(2)}(M,\omega;\alpha_0,\alpha_\infty):=\sup\{\max H\,|
\, H\in {\cal H}_{ad}(M,\omega;\alpha_0,\alpha_\infty)\},\\
C_{HZ}^{(2\circ)}(M,\omega;\alpha_0,\alpha_\infty):=\sup\{\max
H\,| \, H\in{\cal
H}_{ad}^\circ(M,\omega;\alpha_0,\alpha_\infty)\}.
\end{array}
\right.\leqno(1.1)$$
Denote by  $pt$ the generator of $H_0(M)$ represented by a point.
Then we have symplectic invariants
$$
\left\{\begin{array}{ll}
C_{HZ}(M,\omega):=C_{HZ}(M,\omega; pt, pt),\\
C_{HZ}^{\circ}(M,\omega):=\widehat C_{HZ}(M,\omega;pt, pt).
\end{array}
\right.\leqno(1.2)$$ It had been proved in [Lu3] that $C_{HZ}$ is
a symplectic capacity and that ${\cal W}_G\le C_{HZ}\le c_{HZ}$.
In fact, in Lemma 1.4 of the recent [Lu3, v9] we proved that
$C_{HZ}(M,\omega)=c_{HZ}(M,\omega)$ and
$C^\circ_{HZ}(M,\omega)=c^\circ_{HZ}(M,\omega)$ if either $M$ is
closed or each compact subset $K\subset M\setminus\partial M$ may
be contained in a compact submanifold $W\subset M$ with connected
boundary and of codimension zero.  Without special statements we
follow all notations and conventions in [Lu3]. Especially, we
always make the convention that $\sup\emptyset=0$ and
$\inf\emptyset=+\infty$ in this paper. Moreover, we shall omit the
superscripts in $C_{HZ}^{(2)}$ and $C_{HZ}^{(2\circ)}$ without
occurring of confusion.
 \vspace{1mm}

A $2n$-dimensional {\it symplectic toric manifold} is a closed
connected symplectic manifold $(M,\omega)$ equipped with an
effective Hamiltonian action $\tau:\T^n\to {\rm Diff}(M,\omega)$
of the standard (real) $n$-torus $\T^n=\R^n/2\pi\Z^n$ and with a
choice of a corresponding {\it moment map} $\mu:M\to (\R^n)^\ast$.
The image $\triangle=\mu(M)\subset (\R^n)^\ast$ is a convex
polytope, called the {\it moment polytope}. It was proved in [Del]
that the polytope satisfies: (i) there are $n$ edges meeting at
each vertex $p$, (ii) the edges meeting at the vertex $p$ are
rational, i.e., each edge is of the form $p+tv_k$, $0\le
t\le\infty$, where $v_k\in (\Z^n)^\ast$; (iii) the $v_1,\cdots,
v_n$ in $(ii)$ can be chosen to be a basis of $(\Z^n)^\ast$. Such
a polytope is called a {\it Delzant polytope}. It can be uniquely
written as
$$
\triangle=\bigcap^d_{k=1}\Bigl\{ x\in (\R^n)^\ast\,\Bigm|\,
\langle x, u_k\rangle=x(u_k)\ge \lambda_k\Bigl\}.\leqno(1.3)
$$
 Here $d$ is the number of the
$(n-1)$-dimensional faces of $\triangle$, $u_k$ is a uniquely
primitive element of the lattice $\Z^n\subset\R^n$ (the
inward-pointing normal to the $k$-th face of $\triangle$), and
$\lambda_k$ is a real number. Delzant [Del] associated to this
$\triangle\subset (\R^n)^\ast$ a closed connected symplectic
manifold $(M_\triangle,\omega_\triangle)$ of dimension $2n$
together with a Hamiltonian $\T^n$-action
$\tau_\triangle:\T^n\to{\rm Diff}(M_\triangle,\omega_\triangle)$
such that the image of the corresponding moment map
$\mu_\triangle:M_\triangle\to(\R^n)^\ast$ is precisely $\triangle$
and that $(M,\omega,\tau)$ is isomorphic as a Hamiltonian
$\T^n$-space to $(M_\triangle,\omega_\triangle,\tau_\triangle)$.
Two symplectic toric manifolds  are {\it isomorphic} if they are
equivariantly symplectomorphic. Two Delzant polytopes in
$(\R^n)^\ast$ are {\it isomorphic} if they are differ by the
composition of a translation with an element of ${\rm SL}(n,\Z)$.
Delzant showed in [Del] that two symplectic toric manifolds  are
isomorphic if and only if their Delzant polytopes are isomorphic.
Thus we reduce the study of symplectic topology of
$(M,\omega,\tau)$ to that of
$(M_\triangle,\omega_\triangle,\tau_\triangle)$. Delzant's
construction also yielded a ``canonical'' $\T^n$-invariant complex
structure $J_\triangle$ compatible with the symplectic form
$\omega_\triangle$. In other words the quadruple
$(M_\triangle,\omega_\triangle, J_\triangle,\tau_\triangle)$ is a
K\"ahler manifold.

Such manifolds may be explained as a special class of projective
varieties. There exists a class of normal algebraic varieties,
called toric varieties, which are classified by combinatorial
objects called fans.  Let $\Sigma$ be a complete regular fan in
$\R^n$ and $G(\Sigma)=\{u_1,\cdots, u_d\}$ be the set of all
generators of $1$-dimensional cones in $\Sigma$. Denote by
$X_\Sigma$ the compact toric manifold associated with $\Sigma$. If
it is projective, i.e., if $X_\Sigma$ admits a compatible
symplectic structure $\omega$ such that $(X_\Sigma,\omega,J)$ is
K\"ahler, then every K\"ahler form on $X_\Sigma$ can be
represented by a strictly convex support function $\varphi$ for
$\Sigma$ (cf. \S2.2). Conversely every strictly convex support
function for $\Sigma$ represents a K\"ahler form on $X_\Sigma$.
Therefore, in this paper we shall use the same letter to denote a
K\"ahler form on $X_\Sigma$ and the corresponding strictly convex
support function for $\Sigma$ when the context makes our meaning
clear. For such a function $\varphi$ setting
$$\triangle_\varphi=\{x\in (\R^n)^\ast\,|\, \langle x, m\rangle\ge
-\varphi(m)\,\forall m\in\R^n\},\leqno(1.4)
$$
it is a Delzant polytope in $(\R^n)^\ast$. With
$S^1=\{z\in\C\,:\,|z|=1\}$ the action of the maximal compact torus
$\T^n=\R^n/2\pi\Z^n\cong (S^1)^n\subset (\C^\ast)^n$ is
Hamiltonian with respect to the symplectic structure
$2\pi\cdot\varphi$ and has  moment polytope $\triangle_\varphi$.
In other words
$$(X_\Sigma,
2\pi\cdot\varphi)=(M_{\triangle_\varphi},\omega_{\triangle_\varphi}).\leqno(1.5)$$

For the  Delzant polytope $\triangle$ in (1.3) we denote by
$\Sigma_\triangle$ the complete regular fan in $\R^n$ associated
with it and by $P_\triangle=X_{\Sigma_\triangle}$ the
corresponding projective toric manifold (cf. Section 2.2). It
follows from (1.5) that $\Sigma$ is the fan associated with
$\triangle_\varphi$. As showed in [Gu2] the set of all generators
of $1$-dimensional cones in $\Sigma_\triangle$ is given by
$G(\Sigma_\triangle)=\{u_1,\cdots, u_d\}$, and under the identity
$M_\triangle=P_\triangle$ the K\"ahler form $\omega_\triangle$ is
represented by the strictly convex support function for
$\Sigma_\triangle$ defined by $\omega_\triangle(u_i)=-2\pi
\lambda_i$, $i=1,\cdots,d$. Therefore
$\triangle_{\omega_\triangle}=2\pi\triangle$. So we can study the
symplectic topology of toric manifolds from two points of view.
Let $\Z_{\ge 0}$ be the set of all nonnegative integers. Denote
by\vspace{-1mm}
$$\Lambda(\Sigma,\varphi):=\max\sum^d_{i=1}\varphi(u_i)a_i\leqno(1.6)$$\vspace{-1mm}
where $(a_1,\cdots,a_d)\in\Z_{\ge 0}^n$ satisfies
$\sum^d_{i=1}a_iu_i=0$  and $1\le\sum^d_{i=1}a_i\le n+1$. Our
first result is:\vspace{2mm}

\noindent{\bf Theorem 1.1.}\hspace{2mm}{\it  For the above
$\Sigma$ and $\varphi$ one has
$$0<\Lambda(\Sigma,\varphi)\le (n+1)\max_i\varphi(u_i)\quad{\rm
and}\leqno(1.7)$$\vspace{-8mm}
$$ {\cal W}_G(X_\Sigma,\varphi)\le C(X_\Sigma,\varphi; pt,
PD([\varphi]))\le \Lambda(\Sigma,\varphi) \leqno(1.8)
$$
for $C=C_{HZ}^{(2)}$, $C_{HZ}^{(2\circ)}$ and $n\ge 2$. Moreover,
for the interior ${\rm Int}(\triangle_\varphi)$ of
$\triangle_\varphi$ it always holds that
$${\cal W}_G(X_\Sigma,\varphi)\ge \frac{1}{2\pi}{\cal
W}_G({\rm Int}(\triangle_\varphi)\times\T^n, \omega_{\rm
can}),\leqno(1.9)
$$
where
 $({\rm Int}(\triangle_\varphi)\times\T^n,
\omega_{\rm can})=( \{(x,\theta)\,|\, x\in {\rm
Int}(\triangle_\varphi),\;\theta\in\R^n/2\pi\Z^n\},\,
\sum^d_{k=1}dx_k\wedge d\theta_k )$. }\vspace{2mm}

As a by-product of proof of (1.8) we obtain in Corollary 3.2
Mori's theorem on the existence of rational curves through any
point on a uniruled manifold with a different method.
 In Remark 1.5 below we shall give an
example to show that $\Lambda(\Sigma,\varphi)$ may be much smaller
than $ (n+1)\max_i\varphi(u_i)$. In some condition the estimate in
(1.8) can be improved. \vspace{2mm}

\noindent{\bf Theorem 1.2.}\hspace{2mm}{\it If $X_\Sigma$ is also
Fano, i.e., the anticanonical divisor $-K_{X_\Sigma}$ is ample,
then
$$\Upsilon(\Sigma,\varphi):=\inf\Bigl\{\sum^d_{k=1}\varphi(u_k) a_k>0\,
\Bigm|\, \sum^d_{k=1}a_k u_k=0,\,a_k\in\Z_{\ge
0},\,k=1,\cdots,d\Bigr\}>0, \leqno(1.10)
$$
and  for $C=C_{HZ}^{(2)}$, $C_{HZ}^{(2\circ)}$ and any $n\ge 2$,
$$ {\cal W}_G(X_\Sigma,\varphi)\le C(X_\Sigma,\varphi; pt,
PD([\varphi]))\le \Upsilon(\Sigma,\varphi). \leqno(1.11)
$$}\vspace{-4mm}

By the definition it is easy to see that
$\Upsilon(\Sigma,\varphi)\le\Lambda(\Sigma,\varphi)$. In Theorem
2.3 we shall list three equivalent criterions to judge whether or
not  $X_\Sigma$ is Fano from $\Sigma$.

Let $\triangle^n(a):=\{(x_1,\cdots, x_n)\in\R^n_{>0}\,|\,
\sum^n_{k=1}{x_k}<a\}$. For $\triangle$ in (1.3) the following
number
$${\cal W}(\triangle):=\sup\{a>0\,|\, \exists \Psi\in{\rm SL}(n,\Z),
\,x\in(\R^n)^\ast\,s. t.\, \Psi(\triangle^n(a))+x\subset
\triangle\}\leqno(1.12)$$
 is an invariant of the Delzant
polytopes in $(\R^n)^\ast$ under the group generated by elements
of ${\rm SL}(n,\Z)$ and translations. For each vertex $p$ of
$\triangle$ we can assign a positive number $E_p(\triangle)$ as
follows. Let $p_1,\cdots, p_n$ be $n$ vertex adjacent to $p$. By
the above definition of Delzant polytope we may assume that $p_k$
sits in an edge of the form $p+tv_k$, $t\ge 0$, $k=1,\cdots,n$.
Denote by $r_p(\triangle)_k=|p-p_k|/|v_k|$, $k=1,\cdots,n$. Here
$|v|$ denotes the standard norm of vector $v$ in $(\R^n)^\ast$.
Then $p_k=p+r_p(\triangle)_k v_k$, $k=1,\cdots,n$. Let
$$r_p(\triangle)=\bigl\{r_p(\triangle)_1,\cdots,r_p(\triangle)_n\bigr\}
\quad{\rm and}\quad E_p(\triangle)=\min_{1\le k\le
n}r_p(\triangle)_k. \leqno(1.13)
$$
\vspace{2mm}

\noindent{\bf Proposition 1.3.}\quad{\it For the Delzant polytope
$\triangle$ in (1.3) it holds that
$$ \frac{1}{2\pi}{\cal
W}_G({\rm Int}(\triangle)\times\T^n, \omega_{\rm can})\ge {\cal
W}(\triangle)\ge \max_{p\in{\rm
Vert}(\triangle)}E_p(\triangle).\leqno(1.14)$$}

  We also want to derive the estimation in terms of $\triangle$. A
$n$-dimensional integral polytope $\triangle\subset(\R^n)^\ast$
was called {\it reflexive} in [Ba2] if it satisfies: (i) ${\rm
Int}(\triangle)\cap(\Z^n)^\ast=\{0\}$, and (ii) all facets $F$ of
$\triangle$ are supported by an affine hyperplane of the form
$\{m\in(\R^n)^\ast\,|\,\langle m, v_F\rangle =-1\}$ for some
$v_F\in\Z^n$. A equivalent version is that $0\in{\rm
Int}(\triangle)$ and the polar
$\triangle^\circ:=\{x\in\R^n\,|\,\langle m,x\rangle\ge
-1,\,\forall m\in\triangle\}$ is also a $n$-dimensional integral
polytope $\R^n$. A reflexive polytope $\triangle$ is called a {\it
Fano polytope} if the fan $\Sigma_\triangle$ is regular. Clearly,
 a reflexive and Delzant
polytope is also Fano. Note that polytopes $\triangle$ and
$r\cdot\triangle$ yield the same fans for any $r>0$, and that two
toric manifolds corresponding with two isomorphic Delzant
polytopes have  same Fanoness.  Thus a toric manifold
$P_\triangle$  is Fano if and only if $r\cdot (m+\triangle)$ is a
Fano polytope for some $m\in(\R^n)^\ast$ and $r>0$.  In Theorem
2.5 we shall show that the toric manifold $P_\triangle$ associated
with a Delzant polytope $\triangle$ in (1.3) is Fano if and only
if there exist $m\in (\R^n)^\ast$ and $r>0$ such that
$${\rm
Int}(r\cdot(m+\triangle))\cap (\Z^n)^\ast=\{0\}\quad{\rm and}\quad
r\cdot(\lambda_i+\langle m,u_i\rangle)=\pm 1,\;\forall 1\le i\le
d.\leqno(1.15)
$$
 More sufficient and necessary
conditions will be given there. Using this we get the following
corollary of Theorems 1.1 and 1.2. \vspace{2mm}

\noindent{\bf Corollary 1.4}.\quad{\it For the Delzant polytope
$\triangle\subset (\R^n)^\ast$ in (1.3) let $\Lambda(\triangle)$
$(=\Lambda(\Sigma_\triangle,\omega_\triangle ))$ be the maximum of
$-2\pi\sum^d_{i=1}\lambda_ia_i$ for all
$(a_1,\cdots,a_d)\in\Z_{\ge 0}^n$ satisfying
$\sum^d_{i=1}a_iu_i=0$
 and $1\le\sum^d_{i=1}a_i\le n+1$. Then
 $\Lambda(\triangle)\le -2\pi(n+1)\min_i\lambda_i$
 and for $C=C_{HZ}^{(2)}$, $C_{HZ}^{(2\circ)}$ and $n\ge 2$,
$$2\pi{\cal W}(\triangle)\le {\cal
W}_G(M_\triangle,\omega_\triangle)\le
C(M_\triangle,\omega_\triangle; pt, PD([\omega_\triangle]))
\le\Lambda(\triangle). \leqno(1.16)$$ If there exist $r>0$ and
$m\in (\R^n)^\ast$ such that $r\cdot(m+\triangle)$ satisfies
(1.15), then
$$\Upsilon(\triangle):=\inf\Bigl\{-\sum^d_{k=1}\lambda_k
a_k>0\,\Bigm|\, \sum^d_{k=1}a_k u_k=0,\,a_k\in\Z_{\ge
0},\,k=1,\cdots,d\Bigr\}>0, \leqno(1.17)
$$ and for $C=C_{HZ}^{(2)}$, $C_{HZ}^{(2\circ)}$ and any $n\ge 2$ it holds that
$$ {\cal W}_G(M_\triangle,\omega_\triangle)\le
C(M_\triangle,\omega_\triangle; pt, PD([\omega_\triangle]))\le
2\pi\Upsilon(\triangle).
 \leqno(1.18)$$}\vspace{-4mm}

\noindent{\bf Remark 1.5.}\quad The {\it polygon space} associated
with $\alpha=(\alpha_1,\cdots,\alpha_5)=(3/2, 1,1,1, 4/3)$ is a
symplectic toric manifold
 $({\rm Pol}(\alpha),\omega_\alpha)$
 with moment polytope $\triangle_\alpha$ given by
 $$\{(x_1, x_2)\in\R^2\,|\, \frac{1}{2}\le
x_1\le\frac{5}{2},\,
 \frac{1}{3}\le x_2\le\frac{7}{3},\,
 x_1+x_2\ge 1,\,x_1-x_2\ge -1,\, x_2-x_1\ge -1\}.$$
(cf. [HaKn]). Using (1.15) one can prove that it is not Fano. We
can also compute that $\Lambda(\triangle_\alpha)=25\pi/3<15\pi=
-2\pi(n+1)\min_i\lambda_i$. This shows that the second inequality
in (1.7) may be strict. Since $\Psi={\rm diag}(1, -1)\in{\rm
SL}(2,\Z)$ and $\Psi(\triangle^2(1))+(\frac{1}{2}, \frac{3}{2})$
is contained in $\triangle_\alpha$, we get that ${\cal
W}(\triangle_\alpha)\ge 1$. From these we can use (1.16) to obtain
$$2\pi\le {\cal W}_G({\rm
Pol}(\alpha),\omega_\alpha)\le C({\rm Pol}(\alpha),\omega_\alpha;
pt, PD([\omega_\alpha]))\le 25\pi/3.$$ Moreover, we can prove that
$\Upsilon(\triangle_\alpha)=1/6$. So the second inequality in
(1.18), i.e.,
$$
C({\rm
Pol}(\alpha),\omega_\alpha; pt, PD([\omega_\alpha]))\le
2\pi\Upsilon(\triangle_\alpha)=\frac{\pi}{3}
$$
can not hold because the first one in (1.18) always hold. These
show that the second inequalities in (1.11) and (1.18) do not
necessarily hold for non-Fano symplectic toric
manifolds.\vspace{2mm}

Notice that $(a_1,\cdots,a_d)\in\Z_{\ge 0}^n$ is only taken over a
finite set in the definition of $\Lambda(\Sigma,\varphi)$.  Using
the formula in [Sp] it might be possible to get the optimal
estimation for any compact non-Fano toric  manifold.

The rest of the paper is organized as follows. In Section 2 we
give some necessary preliminaries on toric manifolds; the readers
only need to browse through them. The main results are proved in
Section 3. Three examples are given in Section 4.  In Section 5 we
estimate symplectic capacities of the polygon spaces. Finally four
related results are given in Section 6; They are impacts of
symplectic blow-ups on symplectic capacities,  symplectic packings
in toric manifolds and the estimate of Seshadri constants of an
ample line bundle on toric manifolds, and symplectic capacities of
symplectic manifolds with $S^1$-action.\vspace{2mm}

\noindent{\bf Acknowledgement}. I am grateful to Professors V.V.
Batyrev, A. Givental, A. Kresch,  H. Sato, B. Siebert and J.A.
Wi\'sniewski for clarifying some facts. The author also thanks
ICTP at Italy and IHES at Paris  for their financial support and
hospitality.

\section{Preliminaries on toric manifolds}

The basic references for toric manifolds (in alphabetic order) are
[Au], [Ba1], [Ew], [Fu], [Gu2] and [Oda]. The description here
will be presented in the unity notations in [Ba1] and
[Gu2].\vspace{2mm}

\noindent{\bf 2.1. Symplectic toric manifolds.}\quad Let
$(M_\triangle,\omega_\triangle, J_\triangle,\tau_\triangle)$ be
the symplectic toric manifold
 associated with Delzant polytope $\triangle\subset (\R^n)^\ast$ in (1.3), and
 $\mu_\triangle:M_\triangle\to(\R^n)^\ast$ be the moment map of the $\T^n$-action
 $\tau_\triangle$ on it. Denote by $F_k$ the $k$-th $(n-1)$-dimensional
 face of $\triangle$ defined by the equation $\langle x, u_k\rangle=\lambda_k$.
 They yield complex and symplectic submanifolds of $M_\triangle$ of real
codimension $2$,
$$D_1=\mu_{\triangle}^{-1}(F_1),\cdots, D_d=\mu_{\triangle}^{-1}(F_d). \leqno(2.1)$$
  Let $c_k$ be the cohomology class in
$H^2(M_\triangle, \Z)$ dual to $D_k$. The cohomology class
$[\omega_\triangle]$ and the first Chern class of $M_\triangle$
are respectively given by
$$\frac{1}{2\pi}[\omega_\triangle]=-\sum^d_{k=1}\lambda_k c_k\quad{\rm and}\quad
c_1(M_\triangle)=\sum^d_{k=1}c_k\leqno(2.2)
$$
(cf. [Gu1]). As pointed out in [Ab] the arguments in [Gu1] gave a
symplectomorphism
$$({\rm Int}(M_\triangle), \omega_\triangle)\cong
({\rm Int}(\triangle)\times\T^n, \omega_{\rm can}).\leqno(2.3)$$
 Here ${\rm Int}(M_\triangle)=\phi_{\T^n}^{-1}({\rm Int}(\triangle))$ is an
open dense subset in $M_\triangle$, $x\in{\rm Int}(\triangle)$,
$\theta\in\R^n/2\pi\Z^n$ and $\omega_{\rm can}=
\sum^d_{k=1}dx_k\wedge d\theta_k$.
 Thus $(x,\theta)$ may be viewed as symplectic coordinates in
${\rm Int}(M_\triangle)$.\vspace{2mm}

\noindent{\bf 2.2. Fans and toric varieties.}\quad  For an integer
$k\ge 1$, a convex subset $\sigma\subset \R^n$ is called a {\it
regular $k$-dimensional cone} if there exists a $\Z$-basis
$v_1,\cdots,v_k,\cdots,v_n$ of $\Z^n$ such that $\sigma=\R_{\ge
0}v_1+\cdots+\R_{\ge 0}v_k$. Such $v_1,\cdots,v_k\in \Z^n$ are
called the {\it integral generators of} $\sigma$. The origin $0\in
\R^n$ is called the {\it regular zero dimensional cone}. The cones
generated by subsets of the integral generators of $\sigma$ are
called the faces of $\sigma$. A finite system
$\Sigma=\{\sigma_1,\cdots,\sigma_s\}$ of regular cones in $\R^n$
is called a {\it complete regular n-dimensional fan} in $\R^n$ if
(i) any face of each cone $\sigma\in\Sigma$ is also in $\Sigma$;
(ii) the intersection $\sigma_1\cap\sigma_2$ of any two cones
$\sigma_1,\sigma_2\in\Sigma$ is a face of each; (iii)
$\R^n=\sigma_1\cup\cdots\cup\sigma_s$. A toric variety is compact
and nonsingular if and only if its corresponding fan is complete
and regular. We always consider such a fan $\Sigma$ below.  The
set of all $k$-dimensional cones of $\Sigma$ is denoted by
$\Sigma^{(k)}$.
 For every $\sigma\in\Sigma^{(1)}$ there is a
unique generator $u\in\Z^n$ such that $\sigma=\Z_{\ge 0}\cdot u$.
Denote by $G(\Sigma)=\{u_1,\cdots, u_d\}$ the set of all
generators of elements of $\Sigma^{(1)}$. A nonempty subset ${\cal
P}=\{u_{i_1},\cdots, u_{i_k}\}\subset G(\Sigma)$ is called a {\it
primitive collection} if it is not the set of generators of a
$k$-dimensional cone in $\Sigma$, while for each generator
$u_{i_l}\in{\cal P}$ the elements of ${\cal
P}\setminus\{u_{i_l}\}$ generate a $(k-1)$-dimensional cone in
$\Sigma$. Since $\Sigma$ is complete there exists a unique cone
$\sigma({\cal P})\in\Sigma$ whose relative interior contains
$u_{i_1}+\cdots +u_{i_k}$. Let $G(\sigma({\cal
P}))=\{u_{j_1},\cdots, u_{j_m}\}$. We get a linear relation
$$u_{i_1}+\cdots +u_{i_k}=c_{j_1}u_{j_1}+\cdots
+c_{j_m}u_{j_m},\,c_{j_s}>0,\,c_{j_s}\in\Z.\leqno(2.4)$$
 (we allow $m=0$ if $u_{i_1}+\cdots +u_{i_k}=0$.)
 It is called the {\it primitive relation} for ${\cal P}$.
 The integer $\deg({\cal P}):=k-(c_1+\cdots +c_m)$ is called the
 {\it degree} of ${\cal P}$. Denote by ${\rm PC}(\Sigma)$
 the set of primitive collections of $\Sigma$.
Let ${\rm A}({\cal P})=\{z\in\C^d\,|\, z_i=0\;{\rm if}\;
u_i\in{\cal P}\}$ and $Z(\Sigma)=\cup_{{\cal P}}{\rm A}({\cal
P})$, where ${\cal P}$ takes over ${\rm PC}(\Sigma)$. Put
$U(\Sigma)=\C^d\setminus Z(\Sigma)$ and
$$R(\Sigma)=\{\mu=(\mu_1,\cdots,\mu_d)\in\Z^d\,|\, \mu_1 u_1+\cdots
+\mu_d u_d=0\}.$$
 Clearly, $R(\Sigma)$ is isomorphic to
$\Z^{d-n}$. Let $(e_1,\cdots, e_d)$ be the standard basis of
$\R^d$.  Define a linear map $\beta:\R^d\to\R^n,\; e_k\mapsto
u_k,\;\;k=1,\cdots, d$. It maps $\Z^d$ onto $\Z^n$.  Note that the
map $\beta$  can be naturally extended to a map
$\beta_{\C}:\C^d\to\C^n$ that maps $2\pi i\Z^d$ onto $2\pi i\Z^n$.
We still denote by $\beta_{\C}$ the induced map from
$T^d_\C:=\C^d/2\pi i\Z^d$ to $T^n_\C:=\C^n/2\pi i\Z^n$. Let
$N_\C(\Sigma)$ be the kernel of this map. Using the  group
isomorphism ${\rm E}_d: T^d_\C\to (\C^\ast)^d$ given by
$$[w]=[(w_1,\cdots, w_d)]\mapsto (e^{w_1},\cdots, e^{w_d}),$$
we get a subgroup of $(\C^\ast)^d$, ${\rm D}(\Sigma):={\rm
E}_d(N_\C(\Sigma))$. Explicitly,  it is isomorphic to
$(\C^\ast)^{d-n}$ as the Lie group. Moreover,  ${\rm D}(\Sigma)$
acts freely and properly on $U(\Sigma)$. Thus the quotient
$X_\Sigma=U(\Sigma)/{\rm D}(\Sigma)$ is a simply connected compact
complex manifold of dimension $n$, called  the {\it compact toric
manifold associated with $\Sigma$.}  Denote by
$$D_k(\Sigma)=\{[(z_1,\cdots, z_d)]\in U(\Sigma)/{\rm D}(\Sigma)\,|\,
z_k=0\},\;k=1,\cdots,d.\leqno(2.5)$$
 They are complex submanifolds of $X_\Sigma$ of codimension one and
  form a basis for the group ${\rm T_NDiv}(X_\Sigma)$ of
   $T_N=(\C^\ast)^n$-invariant divisors.

A continuous function $\varphi:\R^n\to\R$ is called $\Sigma$-{\it
piecewise linear} if it is a linear function on every cone of
$\Sigma$. Such a function is uniquely determined by its values on
elements $u_k\in G(\Sigma)$. We also call $\varphi\in{\rm
PL}(\Sigma)$ {\it integral} if $\varphi(\Z^n)\subset\Z$. Denote by
${\rm PL}(\Sigma)$ the space of all $\Sigma$-piecewise linear
functions on $\R^n$.  For $\varphi\in {\rm PL}(\Sigma)$ and
$\mu\in R(\Sigma)\otimes\R$ the {\it degree of $\mu$ relative to}
$\varphi$ is defined by
$\deg_\varphi(\mu)=\sum^d_{k=1}\mu_k\varphi(u_k)$
\vspace{2mm}

\noindent{\bf Theorem 2.1.}\hspace{2mm}{\it For $A\in
H_2(X_\Sigma,\Z)$ let $\mu_k(A)$ denote the intersection numbers
$A\cdot D_k(\Sigma)$, $k=1,\cdots, d$. Then
$\mu(A)=(\mu(A)_1,\cdots,\mu(A)_d)\in R(\Sigma)$ and the map
$$H_2(X_\Sigma,\Z)\to R(\Sigma),\; A\mapsto \mu(A)\leqno(2.6)$$
is an isomorphism. Denote by $\Xi_\Sigma$ the inverse of the
isomorphism and its natural extension $R(\Sigma)\otimes\R\to
H_2(X_\Sigma,\R)$. Moreover, the homomorphism $\varphi\mapsto
\sum^d_{k=1}\varphi(u_k)PD(D_k(\Sigma))$ from ${\rm PL}(\Sigma)$
to $H^2(X_\Sigma,\R)$ also induces an isomorphism
$$\Xi^\Sigma: {\rm PL}(\Sigma)/M_{\R}\to H^2(X_\Sigma,\R).\leqno(2.7)$$
In particular, under the isomorphism $\Xi^{\Sigma}$ the first
Chern class $c_1(X_{\Sigma})$ is represented by the class of
$\varphi_{c_1}\in {\rm PL}(\Sigma)$
 such that
$\varphi_{c_1}(u_1)=\cdots=\varphi_{c_1}(u_d)=1$. Furthermore, the
degree-mapping induces the nondegenerate pairing $\deg:{\rm
PL}(\Sigma)/M_{\R}\times R(\Sigma)\otimes\R\to\R$ which coincides
with the canonical intersection pairing $H^2(X_{\Sigma},\R)\times
H_2(X_{\Sigma},\R)\to\R$.}\vspace{2mm}

A nonzero homology class $A\in H_2(X_\Sigma,\Z)$ is called {\it
very effective} in [Kr] if $A\cdot D\ge 0$ for every toric divisor
$D$.  Let ${\rm VNE}(X_\Sigma)$ denote the set of very effective
curve classes on $X_\Sigma$.  Then under the isomorphism (2.6) it
is given by  ${\rm VNE}(X_\Sigma)=\Z^d_{\ge 0}\cap
(R(\Sigma)\setminus\{0\})$.

For each cone  $\sigma=\langle
u_{i_1},\cdots,u_{i_{n-1}}\rangle\in\Sigma^{(n-1)}$ let $\langle
u_{i_1},\cdots,u_{i_{n-1}},u_{i_n}\rangle$ and $\langle
u_{i_1},\cdots,u_{i_{n-1}},u_{i_{n+1}}\rangle$ are the
$n$-dimensional cones in $\Sigma$ which contains $\sigma$ as a
face. Then there are unique integers $b_i\in\Z$, $i=1,\cdots, n+1$
with $b_n=b_{n+1}=1$, such that $b_1u_{i_1}+\cdots+ b_nu_{i_n}+
b_{n+1}u_{i_{n+1}}=0$. We define $v(\sigma)=(v(\sigma)_1,\cdots,
v(\sigma)_d)\in R(\Sigma)$ by $v(\sigma)_r=b_t$ for $r=i_t$ and
$1\le t\le n+1$, and by $v(\sigma)_r=0$ otherwise. Under the
isomorphism (2.6) it corresponds to the class in
$H_2(X_\Sigma,\Z)$ represented by the $T_N$-stable closed
subvariety $V(\sigma)\cong\CP^1$. So the intersection number is
$$(D_l(\Sigma)\cdot V(\sigma))=\left\{\begin{array}{ll}
 b_t & l=i_t\;(1\le t\le n-1)\\
  0& {\rm otherwise}
 \end{array}\right.\leqno(2.8)$$
If $X_\Sigma$ is projective the effective cone  is given by ${\rm
NE}(X_\Sigma)=\sum_{\sigma\in\Sigma^{(n-1)}}\R_{\ge 0}v(\sigma)$.

A $\Sigma$-piecewise linear function $\varphi\in{\rm PL}(\Sigma)$
is called {\it strictly convex support function} for $\Sigma$ if
(i) it is {\it upper convex}, i.e., $\varphi(x)
+\varphi(y)\ge\varphi(x+y)\;\forall x, y\in\R^n$, and (ii)  the
restrictions of it to any two different $n$-dimensional cones
$\sigma_1, \sigma_2\in \Sigma$, are two different linear
functions. Denote by $\varphi_l\in{\rm PL}(\Sigma)$ the unique
functions determined by $\varphi_l(u_k)=\delta_{kl}$,
$k,l=1,\cdots, d$. It is easily checked that they are all upper
convex. Moreover, under the isomorphism (2.7) the divisor
$D_l(\Sigma)\in H^2(X_\Sigma,\R)$ corresponds to the class
represented by $\varphi_l$. Denote by $K(\Sigma)$ the cone in
$H^2(X_\Sigma,\R)\cong{\rm PL}(\Sigma)/(\R^n)^\ast$ consisting of
the classes of all upper convex $\varphi\in{\rm PL}(\Sigma)$, and
by $K^{\circ}(\Sigma)$ the interior of $K(\Sigma)$, i.e., the cone
consisting of the classes of all strictly convex support functions
$\varphi\in{\rm PL}(\Sigma)$. Then $K^{\circ}(\Sigma)\ne\emptyset$
if and only if $X_\Sigma$ is projective.\vspace{2mm}

\noindent{\bf Theorem 2.2.}\quad{\it For a complete regular fan
$\Sigma$ in $\R^n$, $\varphi\in{\rm PL}(\Sigma)$ is a strictly
convex support function for it if and only if the following
equivalent conditions hold.
\begin{description}
\item[(i)] For any primitive collection ${\cal
P}=\{u_{i_1},\cdots, u_{i_k}\}\subset G(\Sigma)$ it holds that
$$\varphi(u_{i_1})+\cdots +\varphi(u_{i_k})>\varphi(u_{i_1}+\cdots +
u_{i_k}).$$
 \item[(ii)] $\triangle_\varphi:=\{m\in (\R^n)^\ast\,|
\,\langle m, x\rangle\ge-\varphi(x),\;\forall x\in \R^n\}$ is a
  Delzant polytope in $(\R^n)^\ast$.
In this case,  for each maximal cone $\sigma\in\Sigma$ let
$\varphi_\sigma\in(\R^n)^\ast$ be the unique element such that
$\langle\varphi_\sigma, x\rangle=-\varphi|_\sigma(x)\,\forall
x\in\sigma$, then different maximal cone give different
$\varphi_\sigma\in(\R^n)^\ast$ and $\{\varphi_\sigma\,|\,
\sigma\in\Sigma^{(n)}\}$ is exactly the set of vertexes of
$\triangle_\varphi$.

\item[(iii)] The divisor $\sum^d_{k=1}\varphi(u_k)D_k$ is ample,
or equivalently
$$((\sum^d_{l=1}\varphi(u_l)D_l(\Sigma))\cdot V(\sigma))=
\sum^d_{k=1}\varphi(u_k)v(\sigma)_k>0\;\; {\rm for}\;{\rm
all}\;\sigma\in\Sigma^{(n-1)}.$$
\end{description}}

(i) is Theorem 4.6 in [Ba1]. The first claim in (ii) follows from
Corollary 2.15 in [Oda], and the second is Lemma 2.12 in [Oda].
(iii) is Theorem 2.18 in [Oda].

With the above fan $\Sigma$ one can associate a polytope
 in $\R^n$
$$\triangle_\Sigma:=\bigcup_{\langle u_1,\cdots,u_k\rangle\in\Sigma}
{\rm conv}(0, u_1,\cdots,u_k)\leqno(2.9)$$
 where $u_i\in G(\Sigma)$ and $\langle u_1,\cdots,u_k\rangle$ is the convex cone
spanned on vectors $u_1,\cdots,u_k$.\vspace{2mm}

\noindent{\bf Theorem 2.3}.\quad{\it The compact toric manifold
$X_\Sigma$ is Fano if and only if the following equivalent
conditions hold.
 \begin{description}
 \item[(i)] $\Sigma$-piecewise linear function $\varphi_{c_1}$
 defined in Theorem 2.1 is strictly convex for $\Sigma$.
 \item[(ii)] Every primitive collection ${\cal P}$ of $\Sigma$ has
 positive degree.
 \item[(iii)] The polytope $\triangle_\Sigma$ is strictly convex
 in the sense that each face of it is of the form ${\rm
 conv}(u_{i_1},\cdots,u_{i_k})$ where $\langle u_{i_1},\cdots,
 u_{i_k}\rangle\in\Sigma$.
 \end{description}}\vspace{2mm}

(i) and (ii) come from [Ba1] and [Ba3] respectively. (iii) was
obtained on page 268 in [Wi].

 There are only finitely many toric
Fano varieties of dimension $n$ up to isomorphism. Toric Fano
manifolds have been classified in low dimensions: there exist
exactly $5$ different toric Del Pezzo surfaces, exactly $18$
different toric Fano $3$-folds and exactly $124$ different toric
Fano $4$-folds. (see [Ba2], [Ba3], [Oda] and references therein).

Since a compact nonsingular toric variety $X_\Sigma$ is projective
(or K\"ahler) if and only if its fan $\Sigma$ comes from some
Delzant polytope, we recall the construction of the fan
$\Sigma_\triangle$ associated with the Delzant polytope
$\triangle$ in (1.3). For each face $F$ of $\triangle$ of
codimension $k$ there exists a unique multi-index $I_F$ of length
$k$, $I_F=(i_1,\cdots, i_k),\quad 1\le i_1<\cdots < i_k\le d$,
such that $F=\{ x\in(\R^n)^\ast\,|\, x(u_i)=\lambda_i,\;\forall
i\in I_F\}$. One has a regular $k$-dimensional cone in $\R^n$,
$\sigma_F=\{\sum t_iu_i\,|\, t_i\ge 0\;\,\forall i\in I_F\}$ with
generators $\{u_i\,|\, i\in I_F\}$. The origin $0\in\R^n$ is
called the regular $0$-dimensional cone. Then the set
$\Sigma_\triangle:=\{\sigma_F\,|\, F\;{\rm is}\;{\rm a}\;{\rm
face}\;{\rm of}\; \triangle\}$ is a complete regular
$n$-dimensional fan in $\R^n$ with
$G(\Sigma_\triangle)=\{u_1,\cdots, u_d\}$, and the corresponding
toric manifold $P_\triangle:=X_{\Sigma_\triangle}$ is projective.
Audin showed in [Au] that there exists a biholomorphism from
$(M_\triangle, J_\triangle)$ to $P_{\triangle}=
U(\Sigma_\triangle)/{\rm D}(\Sigma_\triangle)$ which maps $D_k$ in
(2.1) to $D_k(\Sigma_\triangle)$ in (2.5), $k=1,\cdots, d$. Later
we shall not distinguish between $M_\triangle$ and $P_\triangle$
without special statements.  By Theorem 2.2(ii), if $X_\Sigma$ is
projective then any $\varphi\in K^{\circ}(\Sigma)$ yields a
Delzant polytope $\triangle_\varphi$. It is easily proved that the
fan $\Sigma$ associated with  $\triangle_\varphi$ is exactly
$\Sigma$. Moreover, for any $m\in(\R^n)^\ast$ and $r>0$,  the
above construction implies that
$\Sigma_{m+r\triangle}=\Sigma_\triangle$ and thus
$P_{m+r\triangle}=P_\triangle$ because
$m+r\triangle=\cap^d_{k=1}\{ x\in (\R^n)^\ast\,|\, x(u_k)\ge
m(u_k)+r\lambda_k\}$.\vspace{2mm}

 \noindent{\bf Theorem 2.4.}\hspace{2mm}{\it For the Delzant polytope
 $\triangle$ in (1.3) the following assertions
hold:
\begin{description}
\item[(i)] $K^{\circ}(\Sigma_\triangle)\ne\emptyset$, and the open
cone $K^{\circ}(\Sigma_\triangle)\subset
H^2(P_\triangle,\R)=H^{1,1}(P_\triangle,\R)$ consists of classes
of K\"ahler $(1,1)$-forms on $P_\triangle$. The support function
$h_\triangle:\R^n\to\R$ for $\triangle$ defined by
$$h_\triangle(x)=-\inf\{\langle v, x\rangle\,|\,v\in
\triangle\}\quad\forall x\in \R^n,$$
is strictly convex for
$\Sigma_\triangle$, and $\omega_\triangle=2\pi h_\triangle$.

\item[(ii)] $\triangle_{\omega_\triangle}=2\pi\triangle$, and if
$\varphi\in{\rm PL}(\Sigma_\triangle)$ is strictly convex for
$\Sigma_\triangle$ then one has
$$\triangle_\varphi=\cap^d_{i=1}\{m\in (\R^n)^\ast\,|\,
\langle m, u_i\rangle\ge -\varphi(u_i)\}\quad{\rm and}\quad
(M_{\triangle_\varphi},\omega_{\triangle_\varphi})=(P_\triangle,
2\pi\varphi).$$
\end{description}}\vspace{2mm}

(i) follows from Theorem 2.7 in [Oda]. To prove (ii), it is showed
before that each $(n-1)$-dimensional face
$F_i=\{m\in\triangle\,|\,\langle m, u_i\rangle=\lambda_i\}$ gives
a corresponding $1$-dimensional cone $\sigma_{F_i}=\R_{\ge 0}u_i$
in $\Sigma_\triangle$. By Lemma 2.12 in [Oda] this cone yields a
$(n-1)$-dimensional face
$F^\varphi_i:=\{m\in\triangle_\varphi\,|\,\langle m,
u_i\rangle=-\varphi(u_i)\}$ again.\vspace{2mm}

\noindent{\bf Theorem 2.5.}\quad{\it The projective toric manifold
$M_\triangle=P_\triangle$ is Fano if and only if the following
equivalent conditions hold.
 \begin{description}
\item[(i)] There exist $r>0$ and $m\in (\R^n)^\ast$  such that
${\rm Int}(r\cdot(m+\triangle))\cap (\Z^n)^\ast=\{0\}$ and that
$r\cdot(\lambda_i+\langle m,u_i\rangle)=\pm 1$ for $i=1,\cdots,d$.

\item[(ii)] There  exist $r>0$ and $m\in\R^n$  such that $0\in{\rm
Int}(m+\triangle)$ and that each vertex of $r\cdot(m+\triangle)$
is a primitive vector in $(\Z^n)^\ast$ in the sense that its
coordinates are relatively prime.
\end{description}}\vspace{2mm}

 \noindent{\bf Proof.}\quad   To prove (i), note that
$P_\triangle=P_{\mu(m+\triangle)}$. By Exercise 3.6 on page 70 of
[Gu2], the vertices of the polytope $r\cdot(m+\triangle)$ lie on
integer lattice points if and only if all
$r\cdot(\lambda_i+\langle m, u_i\rangle)$ are integers,
$i=1,\cdots,d$. Moreover it was proved in [Ba2] that
$P_{r\cdot(m+\triangle)}$ is Fano if and only if the integral
polytope $r\cdot(m+\triangle)$ is a reflexive polytope.  These
imply (i). As to (ii) it was proved in [Ew] that
$P_{r\cdot(m+\triangle)}$ is Fano if and only if the integral
polytope $r\cdot(m+\triangle)$ is a Fano ploytope. The condition
in (ii) just right guarantees that $r\cdot(m+\triangle)$ is a Fano
polytope. \hfill$\Box$\vspace{2mm}

\section{Proof of the Main Theorems}

Let $(M,\omega)$ be a closed symplectic manifold.
 For nonzero classes $\alpha_0$, $\alpha_\infty\in
H_\ast(M,\Q)$, using the Gromov-Witten invariant  homomorphism
$\Psi_{A, g, m+2}: H_\ast(\overline{\cal M}_{g, m+2};\Q)\times
H_\ast(M;\Q)^{m+2}\to \Q,$ we defined in [Lu3] a number ${\rm
GW}_g(M,\omega;\alpha_0,\alpha_\infty)$ by the infimum of the
$\omega$-areas $\omega(A)$ of the homology classes $A\in
H_2(M;\Z)$ for which  $\Psi_{A, g,
m+2}(\kappa;\alpha_0,\alpha_\infty,\beta_1,\cdots,\beta_m)\ne 0$
for some homology classes $\beta_1,\cdots,\beta_m\in
H_\ast(M;\Q)$, $\kappa\in  H_\ast(\overline{\cal M}_{g, m+2};\Q)$
and integer $m>0$.
 It was proved in  Theorem 1.10 and Remark 1.11 of [Lu3] that for
 $C=C_{HZ}^{(2)}$, $C_{HZ}^{(2\circ)}$,
$$
C(M,\omega;\alpha_0,\alpha_\infty)\le {\rm
GW}_0(M,\omega;\alpha_0,\alpha_\infty)\quad{\rm and}
\leqno(3.1)
$$\vspace{-7mm}
$$
{\rm GW}_g(M,\omega; pt, PD([\omega]))=\inf\{
{\rm GW}_g(M,\omega; pt,\alpha)\,|\,\alpha\in
H_\ast(M,\Q)\}.\leqno(3.2)
$$
 These are the starting points of proof of our main results.
\vspace{2mm}

\noindent{\bf 3.1. Rational curves on uniruled manifolds.}\quad A
smooth projective variety $X$ over $\C$ is called {\bf uniruled}
if it satisfies the following equivalent conditions:
\begin{description}
\item[(i)] There is a nonempty open subset $U\subset X$ such that
for every $x\in U$ there is a morphism $f:\CP^1\to X$ satisfying
$x\in f(\CP^1)$. \item[(ii)] For every $x\in X$ there is a
morphism $f:\CP^1\to X$ satisfying $x\in f(\CP^1)$.
\end{description}
The following proposition is a key to prove Theorem 1.1. Its proof
was actually contained in Kollar's arguments in [Ko] and
Proposition 7.3 in [Lu3]. For convenience of the readers  we shall
prove it in detail.\vspace{2mm}

\noindent{\bf Proposition 3.1.}\quad{\it Let $X$ be a uniruled
manifold of positive dimension $n$. Then there exist homology
classes $A\in H_2(X;\Z)$ with $1\le c_1(A)\le n+1$, $\alpha\in
H_{2n-2}(X,\Q)$ and $\beta\in H_\ast(X;\Q)$ such that
$$\Psi_{A, 0, 3}(pt; pt,\alpha,\beta)\ne
0.\leqno(3.3)$$}\vspace{-4mm}

\noindent{\bf Proof.}\quad Firstly, note that (3.3) and dimension
condition in the definition of GW- invariants imply
$$2+2n+ (2n-\dim\beta)=2n+ 2c_1(A).$$
It follows that $1\le c_1(A)\le n+1$ because $0\le\dim\beta\le
2n$. So we only need to prove (3.3).

Our proof ideas are based on
 the proof of Theorem 4.2.10 in [Ko] and simple arguments of Gromov-Witten
 invariants.  Recall
the proof of Theorem 4.2.10 in [Ko]. Fix a very general point
$x\in X$ and a very ample divisor $H\subset X$. Since $X$ is
uniruled there exists a rational curve $C$ through $x$ such that
$(C\cdot H)$ is minimal. Let $B=:[C]$. Fix a point $z_0\in\CP^1$
and let $k$ be the complex dimension of the space of morphisms
$f:\CP^1\to X$ such that $f_\ast([\CP^1])=B$ and $f(z_0)=x$. Then
$k\ge 2$ because the isotropic subgroup of automorphism group of
$\CP^1$ at $z_0$ has real dimension $4$. Then for general divisors
$H_1,\cdots, H_k$ linearly equivalent to $H$,
$$\Psi_{B, 0, k+1}(pt; pt, H_1,\cdots, H_k)\ne 0.\leqno(3.4)$$
If $k=2$ then (3.3) holds for $A=B$. If $k=3$ it follows from (6)
in [Mc] that
$$\Psi_{B, 0, 4}(pt; pt, H_1, H_2, H_3)=
\sum_{B=B_1+ B_2}\sum_l\Psi_{B_1, 0, 3}(pt;  pt, H_1, e_l)
\Psi_{B_2, 0, 3}(pt; f_l, H_2, H_2)$$
 where $\{e_l\}_l$ is a basis
for the homology $H_\ast(X;\Q)$ and $\{f_l\}_l$ is the dual basis
with respect to the intersection pairing. This identity implies
that $ \Psi_{B_1, 0, 3}(pt;  pt, H_1, e_l)\ne 0$ for some $l$.
Taking $A=B_1$ one gets (3.3) again. If $k\ge 4$  the composition
law of the GW-invariants gives
\begin{eqnarray*}
\hspace{20mm}& &\Psi_{B, 0, k+1}(pt; H_1,\cdots,H_k)\\
&=&\sum_{B=B_1+ B_2}\sum_{a,b}\Psi_{B_1, 0, 4}(pt; pt, H_1,
H_2,\beta_a) \eta^{ab}\Psi_{B_2, 0, k-1}(pt;\beta_b,
H_3,\cdots,H_k).
\end{eqnarray*}
Here $\{\beta_b\}^L_{b=1}$ is a homogeneous basis of
$H_\ast(X,\Q)$. It follows from (3.4) that
$$\Psi_{B_1, 0, 4}(pt; pt, H_1, H_2,\beta_a)\ne 0$$
for some $B_1\in H_2(X;\Z)$ and $1\le a\le L$. As in case $k=3$ we
can also get (3.3). Clearly we has always $(H^\prime\cdot A)\le
(H^\prime\cdot B)$ for any very ample dvisor $H^\prime$ on $X$.
\hfill$\Box$\vspace{2mm}

\noindent{\bf Corollary 3.2.}\quad{\it For a uniruled manifold $X$
of positive dimension $n$, through any general point of $X$ there
is a rational curve $C$ with $0<(-K_X\cdot C)\le
n+1$.}\vspace{2mm}

This result is not new. It is an easy part of the celebrated
Mori's theorem in [Mor1], [Mor2].  For more general versions of
Corollary 3.2 the reader may refer to [KoMor].\vspace{2mm}

\noindent{\bf 3.2. Proof of Theorem 1.1.}\quad Since $X_\Sigma$ is
uniruled, Proposition 3.1 yields homology classes $A\in
H_2(X_\Sigma;\Z)$ with $1\le c_1(A)\le n+1$, $\alpha\in
H_{2n-2}(X_\Sigma,\Q)$ and $\beta\in H_\ast(X_\Sigma;\Q)$ such
that $\Psi_{A, 0, 3}(pt; pt,\alpha,\beta)\ne 0$. Note that the
Gromov-Witten invariants are deformation invariants. For any
$\varphi\in K^\circ(\Sigma)$ one has
$$\langle[\varphi],A\rangle=\sum^d_{i=1}\varphi(u_i)\mu(A)_i>0.$$
Now $K(\Sigma)$ is the closure of $K^\circ(\Sigma)$ in
$H^2(X_\Sigma,\R)$. Therefore
$\langle[\psi],A\rangle=\sum^d_{i=1}\psi(u_i)\mu(A)_i\ge 0$ for
any $\psi\in K(\Sigma)$. In particular we have
$$\langle[\varphi_l],A\rangle=\sum^d_{i=1}\varphi_l(u_i)\mu(A)_i=
\mu(A)_l\ge 0,\;l=1,\cdots,d.\leqno(3.5)$$
 These show that $A$ is very effective. By Theorem 2.1,
 $c_1(A)=\sum^d_{i=1}\mu(A)_i$. So $1\le\sum^d_{i=1}\mu(A)_i\le
 n+1$. By the definition of $\Lambda(\Sigma,\varphi)$ we have
$$0<\langle[\varphi],A\rangle=\sum^d_{i=1}\varphi(u_i)\mu(A)_i\le\Lambda(\Sigma,
\varphi)$$
 and thus ${\rm GW}_0(M,\omega; pt,\alpha)\le\Lambda(\Sigma,\varphi)$.
 Moreover it is clear that
$$\sum^d_{i=1}\varphi(u_i)\mu_i\le\sum_{\varphi(u_i)>0}\varphi(u_i)\mu_i\le
 (n+1)\max_i\varphi(u_i)$$
for each $\mu\in\Z^n_{\ge 0}$ satisfying $\sum^d_{i=1}\mu_iu_i=0$
and $1\le\sum^d_{i=1}\mu_i\le n+1$.
  By (3.1) and (3.2) we get the desired
(1.8).

The proof of (1.9) is direct. Note that $\triangle_\varphi$  may
be written as
$$\triangle_\varphi=\bigcap^d_{k=1}\{ x\in (\R^n)^\ast\,|\,
\langle x, u_k\rangle=x(u_k)\ge -\varphi(u_k)\}.
$$
By Theorem 2.4(ii) it is a Delzant polytope in $(\R^n)^\ast$, and
$(M_{\triangle_\varphi},\omega_{\triangle_\varphi})=(P_{\triangle_\varphi},
2\pi\varphi)=(X_\Sigma, 2\pi\varphi)$.
 Using (2.3) we can give a symplectic embedding from
 $({\rm Int}(\triangle_\varphi)\times\T^n, \omega_{\rm can})$ to
$(M_{\triangle_\varphi}, \omega_{\triangle_\varphi})$. Then (1.9)
follows from these and the monotonicity of symplectic capacities.
\hfill$\Box$\vspace{2mm}

\noindent{\bf 3.3. Proof of Theorem 1.2.}\quad  For every
$A=\Xi_{\Sigma}(a)\in{\rm VNE}(X_\Sigma)$, by Theorem 9.1 in [Ba1]
the moduli space ${\cal M}(A, X_\Sigma)$ consisting of holomorphic
maps $ f:\CP^1\to X_\Sigma$ with $f_\ast([\CP^1])=A$ is
irreducible and the virtual dimension of it is equal to $n+
c_1(X_\Sigma)(A)= n+ \sum^d_{k=1}a_k.$ Denote by
$m=1+\sum^d_{k=1}a_k$ and by $c_k\in H^2(X_\Sigma, \Z)$ the
Poincare dual of $[D_k(\Sigma)]$, $k=1,\cdots,d$. It was stated in
[Ba1] that
$$c_1^{a_1}\cdots c_d^{a_d}=q^A\leqno(3.6)$$
holds in $QH^\ast(X_\Sigma)$. The author incorrectly admitted it
in [Lu2]. Actually one only can prove (3.6) for all $A\in{\rm
VNE}(X_\Sigma)$ in the toric Fano manifolds. The first proof was
given by Givental in [Giv] (also see [Kr] for an elementary proof
for a class of Fano toric manifolds, and [CiS] for another
different proof for Fano toric manifolds with minimal Chern number
at least two). In terms of GW-invariants (3.6) means
$$\Psi^{X_\Sigma}_{A, 0, m+1}(pt; pt,
\underbrace{D_1(\Sigma),\cdots, D_1(\Sigma)}_{a_1},\cdots,
\underbrace{D_d(\Sigma),\cdots, D_d(\Sigma)}_{a_d})=1.\leqno(3.7)
$$
 Its enumerative interpretation is that for
a given  general point $p_0$ on $X_\Sigma$ and generic distinct
points $z_0, z_{k,i}$, $i=1,\cdots, a_k$ and $k=1,\cdots, d$ on
$\CP^1$ there exists precisely one morphism $f\in {\cal M}(A,
X_\Sigma)$ such that $f(z_0)=p_0$ and $f(z_{k,i})\in D_k(\Sigma)$
for $i=1,\cdots, a_k$ and $k=1,\cdots, d$.
 In particular
 $\varphi(A)>0$. But Theorem 2.1 shows that
$\varphi(A)=\sum^d_{k=1}\omega(u_k) a_k$. Therefore for a given
$a\in\Z^d_{\ge 0}\cap R(\Sigma)$, $\Xi_\Sigma(a)\in {\rm
VNE}(X_\Sigma)$ if and only if $\sum^d_{k=1}\varphi(u_k) a_k>0$.
Now (1.10) can easily follow from this and the Gromov compactness
theorem.  Hence (3.2) and (3.7) give
$${\rm GW}_0(X_\Sigma,\varphi; pt, PD([\varphi]))
\le\varphi(A)=\sum^d_{k=1}\varphi(u_k) a_k$$ for any
$A=\Xi_{\Sigma_\triangle}(a)\in{\rm VNE}(X_\Sigma)$, and thus
$${\rm GW}_0(X_\Sigma,\varphi; pt, PD([\varphi]))
\le \Upsilon(\Sigma,\varphi).
$$
This and (3.1) give (1.11). \hfill$\Box$\vspace{2mm}

\noindent{\bf Remark 3.3.}\quad For a symplectic toric manifold
$({\rm Pol}(\alpha),\omega_\alpha)$ in Remark 1.5 it is easily
seen that (3.6) cannot hold for all $A\in{\rm VNE}({\rm
Pol}(\alpha))$.\vspace{2mm}

\noindent{\bf 3.4. Proof of Proposition 1.3.}\quad Denote by
$\omega_0=\sum^n_{k=1}dx_k\wedge d\theta_k$ and $\omega_{\rm
can}=\sum^n_{k=1}dx_k\wedge d\theta_k$ the standard symplectic
form on $\R^{2n}=\R^n\times\R^n$ and its descending symplectic
form on $\R^n\times\T^n=\R^n\times(\R^n/2\pi\Z^n)$ respectively.
For  $a_k>0, b_k>0$, $k=1,\cdots, n$, we also denote by
$$
E(r_1,\cdots, r_n)=\bigl\{(x_1,y_1,\cdots, x_n,
y_n)\in\R^{2n}\,\bigm|\, \sum^n_{j=1}(x_j^2+ y_j^2)/r_j^2<
1\bigr\},\leqno(3.8)$$\vspace{-5mm}
\begin{eqnarray*}
 \triangle(a_1,\cdots,
a_n)&=&\bigl\{(x_1,\cdots, x_n)\in\R^n_{>0}\,\bigm|\,
\sum^n_{k=1}x_k/a_k<1\bigr\}\subset\R^n,\\
\Box(b_1,\cdots, b_n)&=&\{(\theta_1,\cdots,\theta_n)\in\R^n\,|\;
0<\theta_k<b_k\,\forall 1\le k\le n\}
\end{eqnarray*}
 and abbreviate
$\triangle^n(a):=\triangle(a_1,\cdots, a_n)$ and
$\Box^n(b):=\Box(b_1,\cdots, b_n)$ if $a_1=\cdots =a_n=a$ and
$b_1=\cdots =b_n=b$. The following two lemmas will be also used in
 $\S6$.\vspace{2mm}

\noindent{\bf Lemma 3.4}. ([Sik])\hspace{2mm}{\it Let $U,
V\subset\R^n$ be two connected open sets with $H^1(U)=0$ and
$H^1(V)=0$. For the symplectic submanifolds $U\times\T^n$ and
$V\times\T^n$ of $(\R^n\times\T^n,\omega_{\rm can})$, the
following two statements are equivalent:
\begin{description}
\item[(i)] $(U\times\T^n, \omega_{\rm can})$ and $(V\times\T^n,
\omega_{\rm can})$
 are symplectomorphic;
\item[(ii)] there exists a unimodular matrix $\Phi\in\Z^{n\times
n}$ and a vector $x\in\R^n$ such that $V=\Phi U+ x$.
\end{description}}\vspace{2mm}

\noindent{\bf Lemma 3.5}. ([Sch, Lemma 3.11])\hspace{2mm}{\it Let
$E(c_1,\cdots, c_n)$ be as above. Then for all $\epsilon>0$,
\begin{description}
\item[(i)] $(E(\sqrt{2a_1}-\epsilon,\cdots, \sqrt{2a_n}-\epsilon),
\omega_0)$ embeds symplectically in $(\triangle(a_1,\cdots,
a_n)\times\Box^n(2\pi),\omega_0)$ in such a way that for all
$\alpha\in (0, 1)$, $\alpha E(\sqrt{2a_1}-\epsilon,\cdots,
\sqrt{2a_n}-\epsilon)$ is mapped into
$((\alpha+\epsilon)\triangle(a_1,\cdots, a_n))\times\Box^n(2\pi)$;
\item[(ii)] $(\triangle(a_1-\epsilon,\cdots,
a_n-\epsilon)\times\Box^n(2\pi), \omega_0)$ embeds symplectically
in $E(\sqrt{2a_1},\cdots, \sqrt{2a_n})$ in such a way that for all
$\alpha\in (0, 1)$, $(\alpha\triangle(a_1-\epsilon,\cdots,
a_n-\epsilon))\times\Box^n(2\pi)$ is mapped into
$(\alpha+\epsilon)E(\sqrt{2a_1},\cdots, \sqrt{2a_n})$.
\end{description}}\vspace{2mm}

Now we are in position to prove Proposition 1.3. By Lemma 3.4, if
$\Psi(\triangle^n(a))+x\subset\triangle$ for some $\Psi\in{\rm
SL}(n,\Z)$ and $x\in(\R^n)^\ast$ then there exists a symplectic
embedding from $(\triangle^n(a)\times\T^n, \omega_{\rm can})$ into
$(\triangle\times\T^n, \omega_{\rm can})$. Moreover, Lemma 3.5 can
give a symplectic embedding from
$(B^{2n}(\sqrt{2a}-\epsilon),\omega_0)$ into
$(\triangle^n(a)\times\Box^n(2\pi),\omega_0)\subset
(\triangle^n(a)\times\T^n,\omega_{\rm can})$ for any given small
$\epsilon>0$. The definition of ${\cal W}(\triangle)$ and the
monotonicity of the symplectic capacities yield the first
inequality in (1.14).

In order to prove the second inequality in (1.14) let $p\in{\rm
Vert}(\triangle)$ such that
$$E_p(\triangle)=\max\{E_q(\triangle)\,|\, q\in{\rm Vert}(\triangle)\}.$$
Suppose that $p_1,\cdots, p_n\in{\rm Vert}(\triangle)$ are the
adjacent $n$ vertexes as described above Proposition 1.3. Then
there exists a unique unimodular matrix $A\in {\rm SL}(n,\Z)$ such
that $Ae_k^\ast=v_k$, $k=1,\cdots,n$. So the map
$$\Phi:(\R^n)^\ast\to (\R^n)^\ast,\; x\mapsto Ax-p\leqno(3.9)$$
maps the vertexes $p$ and $p_1,\cdots, p_n$ to the origin and
$r_p(\triangle)_1e_1^\ast, \cdots, r_p(\triangle)_ne_n^\ast$. It
follows that $\Phi$ maps the convex combination ${\rm conv}(p,
p_1,\cdots, p_n)$ onto ${\rm conv}(0,
r_p(\triangle)_1e_1^\ast,\cdots, r_p(\triangle)_ne_n^\ast)$. Since
${\rm conv}(p, p_1,\cdots, p_n)\subset\triangle$,  the inverse map
$\Phi^{-1}$ of $\Phi$ maps  ${\rm conv}(0,
r_p(\triangle)_1e_1^\ast,\cdots, r_p(\triangle)_ne_n^\ast)$ into
$\triangle$. But $\triangle^n(E_p(\triangle))$ is contained in
 ${\rm conv}(0, r_p(\triangle)_1e_1^\ast,\cdots,
 r_p(\triangle)_ne_n^\ast)$. The second inequality in (1.14) is obtained
 immediately.
 \hfill$\Box$\vspace{2mm}

\section{Examples}

\noindent{\bf Example 4.1.}\quad In [CdFKM] Candelas, de la Ossa,
Font, Katz, and Morrison resolved the curve of $\Z_2$
singularities of the weighted projective space
${\CP}^4(1,1,2,2,2)$ to obtain a compact toric manifold
$X=X_\Sigma$. Here the one-dimensional cones in the fan $\Sigma$
are spanned by
$$u_1=-e_1-2e_2-2e_3-2e_4,\; u_2=e_1,\; u_3=e_2,\; u_4=e_3,\; u_5=e_4,\;
u_6=\frac{1}{2}(u_1+u_2)
$$
for the standard basis $e_1, e_2, e_3,e_4$ in $\C^4$. The maximal
cones of $\Sigma$ are as follows:
$$\begin{array}{ll}
\sigma_1=\langle u_1, u_3,u_4,u_5\rangle,\quad\sigma_5=\langle
u_2, u_3,u_4,u_5\rangle\\
\sigma_2=\langle u_1, u_4,u_5,u_6\rangle,\quad\sigma_6=\langle
u_2, u_4,u_5,u_6\rangle\\
\sigma_3=\langle u_1, u_3,u_5,u_6\rangle,\quad\sigma_7=\langle
u_2, u_3,u_5,u_6\rangle\\
\sigma_4=\langle u_1, u_3,u_4,u_6\rangle,\quad\sigma_8=\langle
u_2, u_3,u_4,u_6\rangle.
\end{array}$$
 The only two primitive collections are $\{u_1, u_2\}$ and $\{u_3, u_4, u_5,
u_6\}$.   By Theorem 2.2 a $\Sigma$-piecewise linear function
$\varphi\in{\rm PL}(\Sigma)$ is a strictly convex support function
for $\Sigma$ if and only if
$\varphi(u_3)+\varphi(u_4)+\varphi(u_5)+\varphi(u_6)>0$ and
$\varphi(u_1)+\varphi(u_2)>2\varphi(u_6)$. Note that
$u_1+u_2=2u_6$.  $\varphi_{c_1}\in{\rm PL}(\Sigma)$ is not
 a strictly convex support function for $\Sigma$. By Theorem 2.3, $X$ is not Fano.
 Let $\omega$ be the unique $\Sigma$-piecewise linear function
determined by $\omega(u_1)=1$, $\omega(u_3)=1$ and $\omega(u_i)=0$
for $i=2,4,5,6$. It is easily checked that it is a strictly convex
support function for $\Sigma$. So $X$ is projective and $\omega$
is a K\"ahler symplectic form. Theorem 1.1 yields
$${\cal W}_G(X,\omega)\le C(X,\omega;pt,
PD([\omega]))\le \Lambda(\Sigma,\omega)=1\leqno(4.1)
$$
for $C=C_{HZ}^{(2)}$, $C_{HZ}^{(2\circ)}$.
 Moreover, by Theorem 2.2(ii) we can calculate all vertexes of
 $\triangle_{\omega}$ as follows:
$$\begin{array}{lll}
t_1=(-1,-1,0,0),\quad &t_2=(1,0,0,0),\quad &t_3=(3,1,1,0),\\
t_4=(1,-1,0,1),\quad &t_5=(0,-1,0,0),\quad &t_6=(0,0,0,0),\\
t_7=(0,-1,1,0),\quad &t_8=(0,-1,0,1).
\end{array}$$
Since the matrix
\[ \Phi=\left(\begin{array}{c}
t_2-t_6\\
t_5-t_6\\
t_7-t_6\\
t_8-t_6
\end{array}\right)=\left(\begin{array}{lccr}
1 & 0 & 0& 0\\
0 & -1 & 0 & 0\\
0 & -1& 1& 0\\
0& -1 & 0 & 1\end{array}\right)\] belongs to ${\rm SL}(4,\Z)$ and
maps ${\rm Cl}(\triangle^n(1))={\rm conv}\{0, e_1,e_2,e_3,e_4\}$
onto
$${\rm conv}\{t_6, t_2, t_5, t_7,
t_8\}\subset\triangle_{\omega}.$$
 It follows from Theorem 1.1 and Proposition 1.3 that
 $${\cal W}_G(X,\omega)\ge\frac{1}{2\pi}{\cal
 W}_G(\triangle_{\omega}\times\T^n,\omega_{\rm
 can})\ge{\cal W}(\triangle_{\omega})\ge 1.$$
Combing (4.1) we get
$${\cal W}_G(X,\omega)=\frac{1}{2\pi}{\cal
 W}_G(\triangle_{\omega}\times\T^n,\omega_{\rm
 can})={\cal W}(\triangle_{\omega})=1.$$

\noindent{\bf Example 4.2.}\quad Let $(\CP^n,\omega_{\rm FS})$ be
$n$-dimensional projective space equipped with the Fubini-Study
$\omega_{\rm FS}$. We assume that $\int_{\CP^1}\omega_{\rm
FS}=2\pi$. Then $(\CP^n,\omega_{\rm FS})$ is a $2n$-dimensional
toric manifold and its Delzant polytope has vertices $q_0=0$ and
$q_i=e_i$, $i=1,\cdots,n$. Here $e_1,\cdots, e_n$ are the standard
basis of $\R^n$ and we have identified $(\R^n)^\ast$ with $\R^n$.
Let $p\in \CP^n$ be a fixed point of action of $T^n$ on it
corresponding vertex $q_n$ under the moment map. Since $\CP^n$ is
Fano it easily follows from Corollary 1.4 that for
$C=C_{HZ}^{(2)}$, $C_{HZ}^{(2\circ)}$,
$${\cal W}_G(\CP^n,\omega_{\rm FS})=
C(\CP^n,\omega_{\rm FS}; pt, PD([\omega_{\rm FS}]))=2\pi.
$$
Now take $\tau\in (0,1)$ and consider the $\tau$-blow up of
$(\CP^n,\omega_{\rm FS})$ at  $p$ we get a symplectic toric
manifold $(\widetilde\CP^n_{\tau},\omega_{\tau})$.  By Theorem
1.12 in [Gu2] the vertices of its Delzant polytope
$\triangle_\tau$ are $q_0=0$, $q_n=\delta e_n$, and $q_i=e_i$,
$q_{n+i}=\delta e_n+\delta e_i$, $i=1,\cdots, n-1$.
 Here $\delta=1-\tau$. It is easy to see that
 $$\triangle_\tau=\bigcap^{n+2}_{k=1}\{x\in\R^n\,|\,
 (x,u_k)-\lambda_k\ge 0\}$$
where $u_i=e_i$ and $\lambda_i=0$, $i=1,\cdots,n$, and
$u_{n+1}=-\sum^m_{i=1}e_i$, $u_{n+2}=-e_n$, $\lambda_{n+1}=-1$ and
$\lambda_{n+2}=-\delta$. Note that all Delzant polytopes
$\triangle_\tau$ generate the same fan. All toric manifolds
$\widetilde\CP^n_{\tau}$ are same as complex manifolds. If
$\tau=1/2$ it is easily checked that
$2(\triangle_{\frac{1}{2}}-(\frac{1}{2},\cdots,\frac{1}{2}))$
satisfies (1.15). So $\widetilde\CP^n_{\frac{1}{2}}$ is Fano. In
particular we get that (3.7) holds for
$X_\Sigma=\widetilde\CP^n_{\frac{1}{2}}$. Note that the K\"ahler
forms on $\widetilde\CP^n_{\frac{1}{2}}$ and
$\widetilde\CP^n_{\tau}$ are deformedly equivalent because they
sit in a K\"ahler cone on a complex manifold. Using the fact that
the Gromov-Witten invariants are symplectic deformation invariants
we still obtain (3.7) for $X_\Sigma=\widetilde\CP^n_{\tau}$. As in
the proof of Theorem 1.2 it follows from this and (3.1) that
$$
C(\widetilde\CP^n,\omega_{\tau}; pt, PD([\omega_\tau]))\le
2\pi\Upsilon(\triangle_\tau)=2\pi \delta
$$
for $C=C_{HZ}^{(2)}$, $C_{HZ}^{(2\circ)}$. On the other hand the
first inequality in (1.16) leads to
$$
{\cal W}_G(\widetilde\CP^n,\omega_{\tau})\ge 2\pi(1-\tau)
$$
because $\triangle^n(1-\tau)\subset\triangle$. Hence
$${\cal W}_G(\widetilde\CP^n,\omega_{\tau})=
C(\widetilde\CP^n,\omega_{\tau}; pt,
PD([\omega_\tau]))=2\pi(1-\tau)\leqno(4.2)
$$
for $C=C_{HZ}^{(2)}$, $C_{HZ}^{(2\circ)}$.
 On the other hand it is easily checked that
 $\Lambda(\triangle)\ge \pi (n+1)(1-\tau)$. So (1.18) gives better
 upper bound than (1.16) for $C_{HZ}(\widetilde\CP^n,\omega_{\tau}; pt,
PD([\omega_\tau]))$.\vspace{2mm}

\noindent{\bf Example 4.3.}\quad Consider the following
$4$-dimensional toric Fano manifold $W$ due to Hiroshi Sato [Sa,
Ex.4.7], which was missed in the table of Batyrev [Ba3]. Let $e_1,
e_2, e_3, e_4$ be the standard basis in $\R^4$. Denote by
$u_1=e_1, u_2=e_2, u_3=-e_1-e_2$ and $u_4=e_3, u_5=e_4,
u_6=-e_3-e_4$. Let $W$ be the equivariant blow-ups of
$\CP^2\times\CP^2$ along three $T_N$-invariant $2$-dimensional
irreducible closed subvarieties $\overline{{\rm orb}(\{u_1,
u_4\})}$, $\overline{{\rm orb}(\{u_2, u_5\})}$ and $\overline{{\rm
orb}(\{u_3, u_6\})}$. The set of all generators of $1$-dimensional
cones in its fan $\Sigma$ is $G(\Sigma)=\{u_1,u_2,
u_3,u_4,u_5,u_6,u_7,u_8, u_9\}$, where $u_7=u_1+u_4$,
$u_8=u_2+u_5$ and $u_9=u_3+u_6$. $\Sigma$ has $23$ maximal cones
as follows:
$$\begin{array}{lll}
\sigma_1=\langle u_1, u_2,u_7,u_8\rangle,\quad &\sigma_2=\langle
u_1, u_2,u_6,u_8\rangle,\quad &\sigma_3=\langle u_1, u_2,u_6,u_7\rangle, \\
\sigma_4=\langle u_1, u_3,u_5,u_7\rangle,\quad &\sigma_5=\langle
u_1, u_3,u_5,u_9\rangle,\quad &\sigma_6=\langle u_1, u_3,u_7,u_9\rangle, \\
\sigma_7=\langle u_1, u_5,u_6,u_8\rangle,\quad &\sigma_8=\langle
u_1, u_5,u_6,u_9\rangle,\quad &\sigma_9=\langle u_1, u_5,u_7,u_8\rangle, \\
\sigma_{10}=\langle u_1,u_6,u_7,u_9\rangle,\quad
&\sigma_{11}=\langle
u_2, u_3,u_4,u_9\rangle,\quad &\sigma_{12}=\langle u_2, u_3,u_8,u_9\rangle, \\
\sigma_{13}=\langle u_2,u_3,u_4,u_5\rangle,\quad
&\sigma_{14}=\langle
u_2, u_4,u_7,u_8\rangle,\quad &\sigma_{15}=\langle u_2, u_4,u_6,u_7\rangle, \\
\sigma_{16}=\langle u_2,u_4,u_6,u_9\rangle,\quad
&\sigma_{17}=\langle
u_2, u_6,u_8,u_9\rangle,\quad &\sigma_{18}=\langle u_3, u_4,u_5,u_7\rangle, \\
\sigma_{19}=\langle u_3,u_4,u_7,u_9\rangle,\quad
&\sigma_{20}=\langle
u_3, u_5,u_8,u_9\rangle,\quad &\sigma_{21}=\langle u_4, u_5,u_7,u_8\rangle, \\
\sigma_{22}=\langle u_4,u_6,u_7,u_9\rangle,\quad
&\sigma_{23}=\langle u_5, u_6,u_8,u_9\rangle.
\end{array}$$
Since $W$ is Fano, $\varphi_{c_1}\in{\rm PL}(\Sigma)$ defined by
$\varphi_{c_1}(u_i)=1, i=1,\cdots,9$ gives a symplectic structure
on $W$. It is easily checked that
$$\Upsilon(\Sigma,\varphi_{c_1})=\inf\Bigl\{\sum^9_{i=1}n_i>0\,\Bigm|\,
\sum^9_{i=1}n_iu_i=0, (n_1,\cdots,n_9)\in(\Z_{\ge
0})^9\setminus\{0\}\Bigr\}=1.$$
 By Theorem 1.2 we get that for $C=C_{HZ}^{(2)}$,
 $C_{HZ}^{(2\circ)}$,
$$ {\cal W}_G(W,\varphi_{c_1})\le C(W,\varphi_{c_1}; pt,
PD([\varphi_{c_1}]))\le 1.\leqno(4.3)$$
 By Theorem 2.2(ii) we can calculate all vertexes of
 $\triangle_{\varphi_{c_1}}$ as follows:

$$\begin{array}{lll}
t_1=(1,1,0,0),\quad &t_2=(1,1,-1,0),\quad &t_3=(1,1,0,-1),\\
t_4=(1,-2,0,1),\quad &t_5=(1,-2,-1,1),\quad &t_6=(1,-2,0,0),\\
t_7=(1,0,-2,1),\quad &t_8=(1,-1,-2,1),\quad &t_9=(1,0,0,1),\\
t_{10}=(1,-1,0,-1),\quad &t_{11}=(-2,1,1,-1),\quad &t_{12}=(0,-1,-2,-2),\\
t_{13}=(-2,1,1,1),\quad &t_{14}=(0,1,1,0),\quad &t_{15}=(0,1,1,-2),\\
t_{16}=(-2,-1,0,-1),\quad &t_{17}=(-1,1,-1,0),\quad &t_{18}=(0,-1, 1, 1),\\
t_{19}=(0,-1,-1,-1),\quad &t_{20}=(-1,0,-1,1),\quad &t_{21}=(0,0,1,1),\\
t_{22}=(0,0,1,-2),\quad &t_{23}=(0,0,-2,1).\\
\end{array}$$
Note that the matrix
\[ \Phi=\left(\begin{array}{c}
t_2-t_1\\
t_3-t_1\\
t_9-t_1\\
t_{14}-t_1
\end{array}\right)=\left(\begin{array}{lccr}
0 & 0 & -1& 0\\
0 & 0 & 0 & -1\\
0 & -1& 0& -1\\
-1& 0 & 1 & 0\end{array}\right)\] belongs to ${\rm SL}(4,\Z)$ and
maps ${\rm Cl}(\triangle^n(1))={\rm conv}\{0, e_1,e_2,e_3,e_4\}$
onto
$${\rm conv}\{t_1, t_2, t_3, t_9,
t_{14}\}-t_1\subset\triangle_{\varphi_{c_1}}-t_1.$$
 It follows from Theorem 1.1 and Proposition 1.3 that
 $${\cal W}_G(W,\varphi_{c_1})\ge\frac{1}{2\pi}{\cal
 W}_G(\triangle_{\varphi_{c_1}}\times\T^n,\omega_{\rm
 can})\ge{\cal W}(\triangle_{\varphi_{c_1}})\ge 1.$$
Combing (4.3) we arrive at
$$ {\cal W}_G(W,\varphi_{c_1})=C(W,\varphi_{c_1}; pt,
PD([\varphi_{c_1}]))=1\leqno(4.4)
$$
for $C=C_{HZ}^{(2)}$, $C_{HZ}^{(2\circ)}$.

\section{ Symplectic capacities of polygon
spaces}

Let $\alpha=(\alpha_1,\cdots,\alpha_m)\in\R^m_+$. Following [HaKn]
the {\it polygon space} ${\rm Pol}(\alpha)$, {\it abelian polygon
space} ${\rm APol}(\alpha)$ and {\it upper path space} ${\rm
UP}(\alpha)$  for $\alpha$ are  given by
\begin{eqnarray*}
{\rm
Pol}(\alpha)=\Bigl\{(\rho_1,\cdots,\rho_m)\in(\R^3)^m\,\Bigm|\,
\forall i, |\rho_i|=\alpha_i,\,\sum^m_{i=1}\rho_i=0\Bigr\}\Big/
SO_3,\\
{\rm
APol}(\alpha)=\Bigl\{(\rho_1,\cdots,\rho_m)\in(\R^3)^m\,\Bigm|\,
\forall i,
|\rho_i|=\alpha_i,\,\zeta\bigr(\sum^m_{i=1}\rho_i\bigr)=\alpha_m\Bigr\}\Big/
SO_2,\\
{\rm
UP}(\alpha)=\Bigl\{(\rho_1,\cdots,\rho_{m-1})\in(\R^3)^{m-1}\,\Bigm|\,
\forall i,
|\rho_i|=\alpha_i,\,\zeta\bigl(\sum^{m-1}_{i=1}\rho_i\bigr)\ge\alpha_m\Bigr\}\Big/
\sim
\end{eqnarray*}
respectively, where $SO_3$ acts on $(\R^3)^m$ diagonally,
$\zeta:\R^3\to\R$ is the projection $\zeta(x,y,z)=z$, and
$\rho\sim\rho^\prime$ if $\rho=\rho^\prime$ or if
$\zeta\bigl(\sum^{m-1}_{i=1}\rho_i\bigr)=\alpha_m$ and
$[\rho]=[\rho^\prime]$ in ${\rm APol}(\alpha)$. When $\alpha$ is
{\it generic}, i.e., the equation
$\sum^m_{i=1}\varepsilon_i\alpha_i=0$ has no solution with
$\varepsilon_i=\pm 1$ they are respectively closed symplectic
manifolds of dimensions $2(m-3)$, $2(m-2)$ and $2(m-1)$. Moreover
${\rm Pol}(\alpha)$ is a codimension $2$ symplectic submanifold of
${\rm APol}(\alpha)$, and the latter is a codimension $2$
symplectic submanifold of ${\rm UP}(\alpha)$. In particular ${\rm
APol}(\alpha)$ and ${\rm UP}(\alpha)$  are respectively toric
manifolds with moment polytopes
\begin{eqnarray*}
\Xi_\alpha=\Bigl\{(x_1,\cdots,x_{m-1})\in\prod^{m-1}_{i=1}[-\alpha_i,\alpha_i]\,
\Bigm|\,\sum^{m-1}_{i=1}x_i=\alpha_m\Bigr\},\\
\widehat\Xi_\alpha=\Bigl\{(x_1,\cdots,x_{m-1})\in\prod^{m-1}_{i=1}[-\alpha_i,\alpha_i]\,
\Bigm|\,\sum^{m-1}_{i=1}x_i\ge\alpha_m\Bigr\}.
\end{eqnarray*}
Note that $\widehat\Xi_\alpha$ may viewed as a Delzant polytope.
Indeed, with the standard basis $e_1,\cdots, e_{m-1}$ of
${\R}^{m-1}$ we set $u_i=e_i$, $u_{m-1+i}=-e_i$, $i=1,\cdots, m-1$
and $u_{2m-1}=\sum^{m-1}_{i=1}e_i$. Then with
$\lambda_{m-1+i}=\lambda_i=-\alpha_i$, $i=1,\cdots,m-1$, and
$\lambda_{2m-1}=\alpha_m$ we have
$$\widehat\Xi_\alpha=\bigcap^{2m-1}_{k=1}\bigl\{x\in\R^{m-1}\,\bigm|\, (x,
u_k)-\lambda_k\ge 0\bigr\}.
$$
For nonnegative integers $\mu_k$, $k=1,\cdots,2m-1$ the direct
computation gives rise to
$$
\left\{\begin{array}{ll}
\sum^{2m-1}_{k=1}\mu_ku_k=\sum^{m-1}_{i=1}(\mu_i-\mu_{m-1+i}+
\mu_{2m-1})e_i\vspace{1mm}\\
-\sum^{2m-1}_{k=1}\lambda_k\mu_k=\sum^{m-1}_{i=1}\alpha_i(\mu_i+
\mu_{m-1+i})-\alpha_m \mu_{2m-1}.\end{array}\right.
$$
So $\sum^{2m-1}_{k=1}\mu_ku_k=0\Leftrightarrow \mu_{m-1+i}=\mu_i+
\mu_{2m-1}$, $i=1,\cdots, m-1$, and thus
$1\le\sum^{2m-1}_{i=1}\mu_i\le m\Leftrightarrow 1\le
2\sum^{m-1}_{i=1}\mu_i + m\mu_{2m-1}\le m$. In this case
$$-\sum^{2m-1}_{k=1}\lambda_k\mu_k=2\sum^{m-1}_{i=1}\alpha_i
\mu_i+ (\sum^{m-1}_{i=1}\alpha_i-\alpha_m)\mu_{2m-1}.$$
 Setting $b_i=\mu_i$, $i=1,\cdots,m-1$, and $b_m=\mu_{2m-1}$ we
 get
 $$\Lambda(\widehat\Xi_\alpha)=2\pi\max\Bigl\{2\sum^{m-1}_{i=1}\alpha_i b_i+
 (\sum^{m-1}_{i=1}\alpha_i-\alpha_m)b_m\,\Bigm|\,
 1\le 2\sum^{m-1}_{i=1}b_i+ mb_m\le m,\;b_i\in\Z_{\ge 0}\Bigr\}.\leqno(5.1)$$
By (1.7) it is less than or equal to $2m\pi\max\{\alpha_i\,|\,
1\le i\le  m -1\}$. Moreover, it is easily checked that ${\rm
conv}(0,\delta e_1,\cdots,\delta
e_{m-1})\subset\widehat\Xi_\alpha$ for $\delta=\min\{\alpha_i\,|\,
1\le i\le  m -1\}$. By Theorem 1.1 we have: \vspace{2mm}

\noindent{\bf Proposition 5.1.}\quad{\it Let
$\widehat\omega_\alpha$ denote the symplectic form on ${\rm
UP}(\alpha)$. Then for $C=C_{HZ}^{(2)}$, $C_{HZ}^{(2\circ)}$,
$$2\pi\min\{\alpha_i\,|\,
1\le i\le  m -1\}\le {\cal W}_G({\rm
UP}(\alpha),\widehat\omega_\alpha)\le C({\rm
UP}(\alpha),\widehat\omega_\alpha; pt,
PD([\widehat\omega_\alpha])),$$
and for $m\ge 3$ it holds that
$$ C({\rm
UP}(\alpha),\widehat\omega_\alpha; pt,
PD([\widehat\omega_\alpha]))\le \Lambda(\widehat\Xi_\alpha)\le
 2m\pi\max\{\alpha_i\,|\,1\le i\le  m -1\}.\leqno(5.2)$$}

Since $\Xi_\alpha$ is isomorphic to the Delzant polytope
$$\triangle_\alpha=\Bigl\{(y_1,\cdots,y_{m-2})
\in\prod^{m-2}_{i=1}[-\alpha_i,\alpha_i]\,
\Bigm|\,\alpha_m-\alpha_{m-1}\le\sum^{m-2}_{i=1}y_i\le\alpha_m+
\alpha_{m-1}\Bigr\},$$ as in the proof of Proposition 5.1 we can
derive from Theorem 1.1:\vspace{2mm}

\noindent{\bf Proposition 5.2.}\quad{\it Let $\omega_\alpha$
denote the symplectic form on ${\rm APol}(\alpha)$. Then
$$\Lambda(\triangle_\alpha)=\max\biggl\{2\sum^{m-2}_{i=1}\mu_i\alpha_i+
\mu_{m-1}\cdot\Bigl(\sum^{m-1}_{i=1}\alpha_i-\alpha_m\Bigr)+ \mu_m
\cdot
\Bigl(\alpha_{m-1}+\alpha_m-\sum^{m-2}_{i=1}\alpha_i\Bigr)\biggr\},$$
  where $\mu_i\in\Z_{\ge 0}$, $i=1,\cdots, m$, satisfy
  $$
  1\le 2\sum^{m-2}_{i=1}\mu_i+ (m-1)\mu_{m-1}+(m-3)\mu_m\le m-1.
  $$
Moreover, for $m\ge 4$ and $C=C_{HZ}^{(2)}$, $C_{HZ}^{(2\circ)}$
it holds that
$$\left.\begin{array}{ll}
{\cal W}_G({\rm APol}(\alpha),\omega_\alpha)\le C({\rm
APol}(\alpha),\omega_\alpha; pt, PD([\omega_\alpha]))\le
\Lambda(\triangle_\alpha) \\
\hspace{28mm}\le
2(m-1)\pi\max\{\alpha_1,\cdots,\alpha_{m-2},\alpha_{m-1}+\alpha_m
\}.\end{array}\right.\leqno(5.3)
$$}

By Proposition 1.3 in [HaKn], for generic $\alpha\in\R_+^m$ and
$\delta>\sum^{m-1}_{j=1}\alpha_j$ one has a symplectomorphism
$${\rm APol}(\alpha_1,\cdots,\alpha_m)\cong{\rm
Pol}(\alpha_1,\cdots,\alpha_{m-1},\delta+\alpha_m,\delta).$$
 So it follows from Proposition 5.2 that\vspace{2mm}

\noindent{\bf Proposition 5.3.}\quad{\it Let
$\omega^\prime_\alpha$ denote the symplectic form on
 ${\rm Pol}(\alpha)$. Then for every generic $\alpha\in\R_+^m$ satisfying
$\alpha_{m-1}>\alpha_m>\frac{1}{2}\sum^{m-1}_{j=1}\alpha_j$ the
polygon space ${\rm Pol}(\alpha_1,\cdots,\alpha_m)$ is
symplectomorphic to ${\rm
APol}(\alpha_1,\cdots,\alpha_{m-2},\alpha_{m-1}-\alpha_m)$. So
with $\alpha^\prime=(\alpha_1,\cdots,\alpha_{m-1},
\alpha_{m-1}-\alpha_m)$,
$$
\Lambda(\triangle_{\alpha^\prime})=\max\bigg\{2\sum^{m-3}_{i=1}\mu_i\alpha_i+
\mu_{m-2}\cdot\Bigl(\sum^m_{i=1}\alpha_i-2\alpha_{m-1}\Bigr)+
\mu_{m-1} \cdot \Bigl(2\alpha_{m-1}+
2\alpha_{m-2}-\sum^m_{i=1}\alpha_i\Bigr)\biggr\},
$$
  where $\mu_i\in\Z_{\ge 0}$, $i=1,\cdots, m-1$, satisfy
  $$
  1\le 2\sum^{m-3}_{i=1}\mu_i+ (m-2)\mu_{m-2}+(m-4)\mu_{m-1}\le m-2.
  $$
Moreover, for $m\ge 5$ and $C=C_{HZ}^{(2)}$, $C_{HZ}^{(2\circ)}$
it holds that
$$\left.\begin{array}{ll}
{\cal W}_G({\rm Pol}(\alpha),\omega_{\alpha}^\prime)\le C({\rm
Pol}(\alpha),\omega_{\alpha}^\prime; pt,
PD([\omega_{\alpha}^\prime]))\le
\Lambda(\triangle_{\alpha^\prime}) \\
\hspace{26mm}\le
2(m-2)\pi\max\{\alpha_1,\cdots,\alpha_{m-3},\alpha_{m-2}+\alpha_{m-1}-\alpha_m
\}.\end{array}\right.\leqno(5.4)
$$}

It was shown in \S6 of [HaKn] that for generic
$\alpha\in{\R}^5_{+}$ and $\beta\in{\R}^6_{+}$ both ${\rm
Pol}(\alpha)$ and ${\rm Pol}(\beta)$ are toric manifolds if
$\alpha_1\ne\alpha_2$ and $\alpha_4\ne\alpha_5$, and if
$\beta_1\ne\beta_2$ and $\beta_5\ne\beta_6$.

\section{Related results}

\noindent{\bf 6.1. Impacts on symplectic capacities of symplectic
blow-ups.}\quad If we always require the class $\kappa$ to be a
single point one $pt$ in the definition of the number ${\rm
GW}_g(M,\omega;\alpha_0,\alpha_\infty)$ at the beginning of
Section 3, the corresponding infimum is denoted by
$$\widehat{\rm GW}_g(M,\omega;\alpha_0,\alpha_\infty).\leqno(6.1)$$
Then ${\rm GW}_g(M,\omega;\alpha_0,\alpha_\infty)\le
 \widehat {\rm GW}_g(M,\omega;\alpha_0,\alpha_\infty)$.

Firstly, it is easy to see  that the symplectic blow up operation
for a symplectic manifold must decrease the volume of it. One
easily find a noncompact symplectic manifold for which a suitable
symplectic blowing up does not decrease its Gromov symplectic
width. Therefore it is a complicated problem. For simplicity we
restrict our attention to the case of a symplectic blow up of a
closed $2n$-dimensional symplectic manifold $(M,\omega)$  at $k$
distinct points. Let $\psi=\coprod_i\psi_i:
\coprod_i(B^{2n}(r_i),\omega_0)\to (M,\omega)$ be a symplectic
embedding  of $k$ disjoint standard symplectic balls of radii
$r_1,\cdots,r_k$, and $\Theta:(\widetilde
M_\psi,\widetilde\omega_\psi)\to (M,\omega)$ be the symplectic
blow-up associated with $\psi$ at $p_i=\psi_i(0)$, $i=1,\cdots,k$.
Let $H_j(M)$ (resp. $H^j(M)$) denote $H_j(M,\Z)$ (resp.
$H^j(M,\Z)$) modulo torsion. Denote by
$\Sigma_i=\Theta^{-1}(0)\thickapprox\CP^{n-1}$ the exceptional
divisor corresponding to $p_i$. Let $E_1,\cdots, E_k$ denote the
homology classes of the exceptional divisors in
$H_{2n-2}(\widetilde M)$ and $e_1,\cdots, e_k\in H^2(\widetilde
M)$ be their Poincar\'e duals. Then
$$[\widetilde\omega]=[\Theta^\ast\omega]-\sum^k_{i=1}\pi
r_i^2 e_i\leqno(6.2)
$$
in $H^2(\widetilde M,\R)$ ([McP]). Let $E_i^\prime\in
H_2(\widetilde M)$ be the classes of lines in the exceptional
divisors $\Sigma_i$ such that $PD(E_i^\prime)=-(-e_i)^{n-1}$. Let
$\{T_0,\cdots,T_q\}$ be a homogeneous basis of $H^2(M)$ of
increasing codimension such that $T_0$ is the fundamental class
and $T_q=pt$. With $p=q+k(n-1)$ we define $\tilde T_{q+1},\cdots,
\tilde T_p$ to be the classes $e_i^j\in H^2(\widetilde M)$,
$i=1,\cdots,k$ and $j=1,\cdots,n-1$. Denote by $\tilde
T_i=\Theta^\ast T_i$, $i=1,\cdots,q$. Then $\{\tilde T_1,\cdots,
\tilde T_p\}$ is a homogeneous basis of $H^2(\widetilde M)$. The
classes $\tilde T_1,\cdots, \tilde T_q$ (resp. $\tilde
T_{q+1},\cdots,\tilde T_p$) are called {\it non-exceptional}
(resp. {\it exceptional}). Note that
$$\left.\begin{array}{ll}
\tilde T_j\cdot \tilde T_{j^\prime}=\Theta^\ast(T_j\cdot
T_{j^\prime}),\quad \tilde T_j\cdot e_i^m=0\\
e_i^m\cdot
e_{i^\prime}^{m^\prime}=\delta^i_{i^\prime}e_i^{m+m^\prime},\quad
e_i^n=(-1)^{n-1}pt
\end{array}\right.\leqno(6.3)$$
 on $\widetilde M$ for $1\le j,j^\prime\le q$, $1\le
i,i^\prime\le k$ and $1\le m,m^\prime\le n-1$. One  has a
canonical decomposition
$$H_2(\widetilde M)=H_2(M)\oplus\Z E_1^\prime\oplus\cdots\oplus\Z
E_k^\prime.
$$
By (6.3) the classes $PD(E_i^\prime)\in H^2(\widetilde M)$ satisfy
$PD(E_i^\prime)\cdot E_j=E_j(E_i^\prime)=-\delta^i_j$. So one has
$$H^2(\widetilde M,\R)=H^2(M,\R)\oplus\R e_1^{n-1}\oplus\cdots\oplus\R
e_k^{n-1}.
$$
As usual we denote by $\Theta!$ the transfer map $PD_{\widetilde
M}\circ\Theta^\ast\circ PD_M$ from $H_\ast(M)$ to
$H_\ast(\widetilde M)$ and call the image $\Theta!(A)$ the
corresponding non-exceptional class  of $A\in H_2(M)$. Using (6.2)
and (6.3) it is not hard to derive that
$$\langle[\widetilde\omega],\Theta!(A)\rangle=
\langle[\omega], A\rangle.\leqno(6.4)$$ Let $pt_M$ (resp.
$pt_{\widetilde M}$) denote the single point class in $H_0(M)$
(resp. $H_0(\widetilde M)$) such that $\langle
PD_M(pt_M),[M]\rangle=1$ (resp. $\langle PD_{\widetilde
M}(pt_{\widetilde M}),[\widetilde M]\rangle=1$). Note that
$H^{2n}(\widetilde M,\R)=\R [\widetilde\omega^n]$ and
 $H^{2n}(M,\R)=\R [\omega^n]$.
 It is
easily checked that
$$\Theta!(pt_M)=pt_{\widetilde M}\quad{\rm and}\quad
\sum^k_{i=1} r_i^{2n}(-\pi
e_i)^n=-\Bigl(\pi^n\sum^k_{i=1}r_i^{2n})PD_{\widetilde
M}(pt_{\widetilde M}\Bigr).\leqno(6.5)
$$
It was proved in [Ga] and
[Hu] that
$$\Psi^{\widetilde M}_{\Theta!(A), 0, m}(pt;\Theta!(\gamma_1),
\cdots,\Theta!(\gamma_m))=\Psi^M_{A, 0, m}(pt;\gamma_1,
\cdots,\gamma_m)\leqno(6.6)
$$
for any $A\in H_2(M)$ and $\gamma_i\in H_\ast(\widetilde M)$,
$i=1,\cdots,m$. We here use the homology classes for convenience.
It follows from this, (6.1) and (6.4) that
$$\widehat{\rm GW}_0(\widetilde M,\widetilde\omega;
\Theta!(\alpha_0),\Theta!(\alpha_\infty))\le \widehat{\rm
GW}_0(M,\omega;\alpha_0,\alpha_\infty).\leqno(6.7)$$
 Note that the first identity in (6.5) and (6.6) give
$$\Psi^{\widetilde M}_{\Theta!(A), 0, m+1}(pt; pt_{\widetilde M},
\Theta!(\gamma_1), \cdots,\Theta!(\gamma_m)) =\Psi^M_{A, 0,
m+1}(pt; pt_M, \gamma_1, \cdots,\gamma_m).\leqno(6.8)$$
 By (3.2),
 ${\rm GW}_0(\widetilde M,\widetilde\omega; pt_{\widetilde M},
PD([\widetilde \omega]))$ is equal to the infimum of the
$\widetilde\omega$-areas $\widetilde\omega(A)$ of the classes
$A\in H_2(\widetilde M)$ for which  $\Psi^{\widetilde M}_{A, 0,
m+1}(\kappa; pt_{\widetilde M},\beta_1,\cdots,\beta_m)\ne 0$ for
some classes $\beta_1,\cdots,\beta_m\in H_\ast(\widetilde M;\Q)$,
$\kappa\in H_\ast(\overline{\cal M}_{0, m+1};\Q)$ and integer
$m>1$.
 Hence it follows from (6.8) that
$$\left.\begin{array}{ll}
{\rm GW}_0(\widetilde M,\widetilde\omega; pt_{\widetilde M},
PD([\widetilde \omega]))\\
\le\inf\{\omega(A)\,|\, \Psi^M_{A, 0, m+1}(pt; pt_M, \gamma_1,
\cdots,\gamma_m)\ne 0\},\end{array}\right.\leqno(6.9)
$$
where the infimum is taken for $A\in H_2(M)$ and $\gamma_i\in
H_\ast(M)$. By (3.1) and (6.7) we obtain:\vspace{2mm}

\noindent{\bf Theorem 6.1.}\quad{\it For any nonzero classes
$\alpha_0,\alpha_\infty\in H_\ast(M,\Q)$ it holds that
$$\left.\begin{array}{ll}
C_{HZ}(\widetilde M,\widetilde\omega;
\Theta!(\alpha_0),\Theta!(\alpha_\infty))\le \widehat{\rm
GW}_0(M,\omega;\alpha_0,\alpha_\infty)\quad{\rm
and}\\
C_{HZ}(\widetilde M,\widetilde\omega; pt_{\widetilde M},
PD([\widetilde \omega]))\le\inf\{\omega(A)\,|\, \Psi^M_{A, 0,
m+1}(pt; pt_M, \gamma_1, \cdots,\gamma_m)\ne
0\}.\end{array}\right.$$}

Notice that the blow-ups of a toric manifold at its toric fixed
points are also toric manifolds. However, the blow up of a toric
Fano manifold is not necessarily Fano again. By (3.7) and Theorem
6.1 we get:\vspace{2mm}

\noindent{\bf Theorem 6.2.}\quad{\it Let $X_{\widetilde\Sigma}$ be
a toric manifold obtained by a sequence of blowings up of a toric
Fano manifold at toric fixed points. So
$G(\Sigma)=\{u_1,\cdots,u_d\}\subset G(\widetilde\Sigma)$. Then
for any strictly convex support function $\varphi$ for
$\widetilde\Sigma$ (also strictly convex for $\Sigma$) it holds
that
$$ {\cal W}_G(X_{\widetilde\Sigma},\varphi)\le C(X_{\widetilde\Sigma},\varphi; pt,
PD([\varphi]))\le 2\pi\cdot\Upsilon(\Sigma,\varphi). \leqno(6.10)
$$
for every $n\ge 2$ and $C=C_{HZ}^{(2)}$, $C_{HZ}^{(2\circ)}$. Here
$\Upsilon(\Sigma,\varphi)$ is given by (1.10) and is always more
than zero though $\Upsilon(\widetilde\Sigma,\varphi)$ might equal
to zero in the case $X_{\widetilde\Sigma}$ is not
Fano.}\vspace{2mm}

The fan $\widetilde\Sigma$ may be obtained from $\Sigma$ by a
sequence of regular stellar operations. For the Delzant polytope
$\triangle$ in (1.3) and a vertex $p$ of it, let the rays $p+tv_i,
t\ge 0$, form the edges of $\triangle$ at $p$ as in the definition
of Delzant polytope above, we choose
$0<\varepsilon<E_p(\triangle)$ and replace the vertex $p$ by the
$n$ vertices $p+\varepsilon v_i,\;i=1,\cdots,n$ to get a new
Delzant polytope $\triangle_\varepsilon$. Then the  symplectic
toric manifold
$(M_{\triangle_\varepsilon},\omega_{\triangle_\varepsilon})$ is
the symplectic $\varepsilon$ blow-up of the symplectic toric
manifold
$(M_\triangle,\omega_\triangle,\tau_\triangle,\mu_\triangle)$ at a
fixed point $q=\mu_\triangle(p)$  of the $\T^n$-action
$\tau_\triangle$.\vspace{2mm}

\noindent{\bf Corollary 6.3.}\quad{\it Suppose that there exist
$r>0$ and $m\in (\R^n)^\ast$ such that $r\cdot(m+\triangle)$
satisfies (1.15), i.e., $M_\triangle$ is Fano. Then
 for $C=C_{HZ}^{(2)}$, $C_{HZ}^{(2\circ)}$ and any $n\ge 2$ it holds that
$$ {\cal W}_G(M_{\triangle_\varepsilon},\omega_{\triangle_\varepsilon})\le
C(M_{\triangle_\varepsilon},\omega_{\triangle_\varepsilon}; pt,
PD([\omega_{\triangle_\varepsilon}]))\le 2\pi\Upsilon(\triangle).
$$}

\noindent{\bf 6.2. Symplectic packings in symplectic toric
manifolds.}\quad We here presents a symplectic packing result in
symplectic toric manifolds via symplectic ellipsoid of form (3.8)
(see [Bi], [Gr], [McP], [Ka], [Tr], [Sch] and references therein
for the exposition and related results). Denote by ${\rm
Vert}(\triangle)$ the vertex set of the Delzant polytope
$\triangle$ in (1.3). For each $p\in{\rm Vert}(\triangle)$ let
$p_1,\cdots, p_n$ be its adjacent $n$ vertexes.  If $\sharp{\rm
Vert}(\triangle)>n+1$ there must exist another $p^\prime\in{\rm
Vert}(\triangle)$ and adjacent $n$ vertexes $p^\prime_1,\cdots,
p^\prime_n$ corresponding to it such that
$$
({\rm conv}(p, p_1,\cdots, p_n))^\circ\cap ({\rm conv}(p^\prime,
p^\prime_1,\cdots, p^\prime_n))^\circ= \emptyset. \leqno(6.11)$$
Hereafter $S^\circ$ denotes the interior of the set $S$. In this
case we say that the vertexes $q$ and $q^\prime$ are {\it
simplicially separating in $\triangle$}. Notice also that each
${\rm conv}(p, p_1,\cdots, p_n)$ determines a family of open
symplectic ellipsoids
$E(\triangle,p,\epsilon):=E(\sqrt{2r_p(\triangle)_1}-\epsilon,\cdots,
\sqrt{2r_p(\triangle)_n}-\epsilon)$ for $0\le\epsilon\le
E_p(\triangle)$.\vspace{2mm}

\noindent{\bf Theorem 6.4.}\hspace{2mm}{\it If any two points of a
given subset $\{q_1,\cdots, q_m\}\subset{\rm Vert}(\triangle)$ are
simplicially separating in $\triangle$, then for any small
$\epsilon>0$ there exists a symplectic packing of $(M_\triangle,
\omega_\triangle)$ via the ellipsoids $E(\triangle,
q_k,\epsilon)$, $k=1,\cdots,m$.}\vspace{2mm}

\noindent{\bf Proof.}\quad By (2.3) it suffices to prove that for
any small $\epsilon>0$ there exists a symplectic embedding of a
disjoint union of the open ellipsoids $E(\triangle,
q_k,\epsilon)$, $k=1,\cdots,m$.  For each $k=1,\cdots,m$, as in
proof of Proposition 1.3  we have the unimodular matrixes
$A_k\in{\rm SL}(n,\Z)$ such that the corresponding transformations
 $$\Phi_k: (\R^n)^\ast\to (\R^n)^\ast,\; x\mapsto A_kx-q_k,$$
  map $q_k, p_{k1},
\cdots, p_{kn}$  to  $0, a_{k1}e_1^\ast,\cdots, a_{kn}e_n^\ast$,
   $k=1,\cdots, m$,  respectively. Here
$p_{k1},\cdots,p_{kn}$ are the adjacent $n$ vertexes to $q_k$, and
$a_{ki}=r_{q_k}(\triangle)_i$, $i=1,\cdots, n$, and $k=1,\cdots,
m$. Now each $\Phi_k$ induces a symplectomorphism ${\cal A}_k$ of
$((\R^n)^\ast\times\T^n,\omega_{\rm can})$ to itself that maps
${\rm conv}(q_k, p_{k1},\cdots, p_{kn})\times\T^n$ onto ${\rm
conv}(0, a_{k1}e_1^\ast,\cdots, a_{kn}e_n^\ast)\times\T^n$,
$k=1,\cdots, m$. Note that
$$({\rm conv}(0, a_{k1}e_1^\ast,\cdots, a_{kn}e_n^\ast))^\circ=
\triangle(a_{k1},\cdots, a_{kn}),\;\;k=1,\cdots, m, $$ provided
that  $(\R^n)^\ast$ is identified with $\R^n$ by the isomorphism
$$x_1e_1^\ast+\cdots + x_ne_n^\ast\mapsto x_1e_1+\cdots + x_ne_n.$$
We can use Lemma 3.5 to find  symplectic embeddings ${\cal B}_k$
of
$$E(\triangle,
q_k,\epsilon)=E(\sqrt{2a_{k1}}-\epsilon, \cdots,
\sqrt{2a_{kn}}-\epsilon)$$ into $(\triangle(a_{k1}, \cdots,
a_{kn})\times \Box^n(2\pi), \omega_0)$ and thus into
$(\triangle(a_{k1},\cdots, a_{kn})\times \T^n, \omega_{\rm can})$,
$k=1,\cdots, m$. Then it is easily checked that the compositions
${\cal A}_k^{-1}\circ{\cal B}_k$, $k=1,\cdots, m$, give the
desired symplectic embeddings. \hfill$\Box$\vspace{2mm}

\noindent{\bf Remark 6.5.}\hspace{2mm}Let
$\underline{a}=(a_1,\cdots, a_n)$ be a vector of positive weights
and $\triangle^n(\underline{a}):=\triangle(a_1,\cdots, a_n)$. Also
denote by $E(\sqrt{2\underline{a}}):=E(\sqrt{2a_1},\cdots,
\sqrt{2a_n})$. The above proof actually shows that if for some
$\triangle^n(\underline{a}^{(k)})\subset \R^n\equiv (\R^n)^\ast$
there exist $A_k\in{\rm SL}(n,\Z)$ and $q_k\in (\R^n)^\ast$,
$k=1,\cdots, m$, such that the sets
$A_k(\triangle^n(\underline{a}^{(k)}))+ q_k\subset \triangle$,
$k=1,\cdots, m$, are mutually disjoint, then $(M_\triangle,
\omega_\triangle)$ admits a symplectic packing via $m$ open
ellipsoids $E(\sqrt{2\underline{a}^{(k)}})$, $k=1,\cdots,
m$.\vspace{2mm}

\noindent{\bf Example 6.6.} Consider the polygon space $({\rm
Pol}(\alpha),\omega_\alpha)$ in Remark 1.5. Its moment polytope
$\triangle_\alpha$ has vertexes: $q_1=(1/2, 3/2)$, $q_2=(4/3,
7/3)$, $q_3=(5/2, 7/3)$, $q_4=(5/2, 3/2)$, $q_5=(4/3, 1/3)$,
$q_6=(2/3, 1/3)$ and $q_7=(1/2,1/2)$. It is easily computed that
$$\left.\begin{array}{ll}
E(\triangle,q_1)=E(\sqrt{2},\sqrt{5\sqrt{2}/3}),\quad
 E(\triangle, q_2)=E(\sqrt{5\sqrt{2}/3}, \sqrt{7/3}),\\
E(\triangle,q_3)=E(\sqrt{7/3},\sqrt{5/3}),\quad
 E(\triangle, q_4)=E(\sqrt{5/3}, \sqrt{7\sqrt{2}/3}),\\
E(\triangle,q_5)=E(\sqrt{7\sqrt{2}/3},\sqrt{4/3}),\quad
 E(\triangle, q_6)=E(\sqrt{4/3}, \sqrt{\sqrt{10}/3}),\\
 E(\triangle,q_7)=E(\sqrt{\sqrt{10}/3},\sqrt{2}).
 \end{array}\right.$$
By Theorem 6.4, for any $\epsilon>0$ sufficiently small, $({\rm
Pol}(\alpha),\omega_\alpha)$ admits the symplectic packings via
the following groups of ellipsoids:
$$\left.\begin{array}{ll}
\{E(\triangle,q_1,\epsilon),E(\triangle,q_3,\epsilon),E(\triangle,q_5,\epsilon)\},\quad
\{E(\triangle,q_1,\epsilon),E(\triangle,q_3,\epsilon),E(\triangle,q_6,\epsilon)\},\\
\{E(\triangle,q_2,\epsilon),E(\triangle,q_4,\epsilon),E(\triangle,q_6,\epsilon)\},\quad
\{E(\triangle,q_2,\epsilon),E(\triangle,q_4,\epsilon),E(\triangle,q_7,\epsilon)\},\\
\{E(\triangle,q_3,\epsilon),E(\triangle,q_5,\epsilon),E(\triangle,q_7,\epsilon)\},\quad
\{E(\triangle,q_1,\epsilon),E(\triangle,q_4,\epsilon),E(\triangle,q_6,\epsilon)\},\\
\{E(\triangle,q_2,\epsilon),E(\triangle,q_5,\epsilon),E(\triangle,q_7,\epsilon)\}.
\end{array}\right.$$

 \noindent{\bf 6.3. Seshadri constants.}\quad
 Let $(M, J)$ be a compact complex manifold of dimension $n$, and $L\to M$
an ample line bundle. Demailly [Dem] defined the  Seshadri
constant of $L$ at a point $x\in M$ to be the infimum
$\varepsilon(L, x)$ of $\int_C c_1(L)/{\rm mult}_xC$,
   where $C$  takes over all irreducible curves passing
through the point $x$, and ${\rm mult}_xC$ is the multiplicity of
$C$ at $x$.   The global {\it Seshadri constant} is defined by
$\varepsilon(L):=\inf_{x\in M}\varepsilon(L, x)$.
 For the toric manifold $X_\Sigma$ as in Theorem 1.1 let $L_k\to P_\triangle$  be
the corresponding line bundles to the toric divisors $D_k(\Sigma)$
in (2.5), $k=1,\cdots, d$. It is well-known that the Chern class
$c_1(L_k)$ is Poincar\'e dual to $[D_k]\in H_2(X_\Sigma,\Z)$ for
each $k$. \vspace{2mm}

\noindent{\bf Theorem 6.7}\hspace{2mm}{\it Let $\Sigma$ be a
complete regular fan in $\R^n$. Then for any ample line bundle
$L\to X_\Sigma$ and any strictly convex support function
 $\varphi_L$  representing the class $c_1(L)$ it holds that
$$\varepsilon(L)\le 2\pi\cdot\Lambda(\Sigma,\varphi_L).\leqno(6.12)$$
Furthermore, if $X_\Sigma$ is also Fano then
$$\varepsilon(L)\le
2\pi\cdot\Upsilon(\Sigma,\varphi_L).\leqno(6.13)$$}

\noindent{\bf Proof.}\hspace{2mm} Recall that in Definition 1.26
of [Lu3, v9] we defined
$${\rm GW}(M,\omega)=\inf{\rm GW}_g(M,\omega; pt,\alpha),$$
where the infimum is taken over all nonnegative integers $g$ and
all homology classes $\alpha\in H_\ast(M;\Q)\setminus\{0\}$ of
degree $\deg\alpha\le\dim M-1$. Using Proposition 6.3 in [BiCi] we
showed in Theorem 1.36 of [Lu3, v9] that for a closed connected
complex manifold $(M, J)$ of dimension $\dim_{\R}M>2$ and any
ample line bundle $L\to M$ it holds that $\varepsilon(L)\le {\rm
GW}(M,\omega_L)$. Here $\omega_L$ is any $J$-compatible K\"ahler
form (the curvature form for a suitable metric connection on $L$)
representing the cohomology class $c_1(L)$. From the proofs of
Theorems 1.1 and 1.2 it is easily seen that for a toric manifold
$X_\Sigma$ and a strictly convex support $\varphi$ for $\Sigma$
one has
$${\rm GW}(X_\Sigma,\varphi)\le 2\pi\Lambda(\Sigma,\varphi)\quad{\rm and}\quad
{\rm GW}(X_\Sigma,\varphi)\le 2\pi\Upsilon(\Sigma,\varphi)$$ in
general case and Fano case respectively. \hfill$\Box$\vspace{2mm}

\noindent{\bf 6.4. Symplectic capacities of symplectic manifolds
with $S^1$-action.}\quad The symplectic toric manifolds are a
special class of symplectic manifolds with the Hamiltonian
$S^1$-action. Let $\{\lambda_t\}=\lambda: S^1=\R/\Z\to {\rm
Ham}(M,\omega)$ be a Hamiltonian circle action on a connected
symplectic manifold of dimension $2n$. Let $H:M\to\R$ be the
Hamiltonian function for the action. It means that the circle
action is generated by the Hamiltonian vector field $X_H$. This
action is called semi-free if it is free on $M\setminus M^{S^1}$.
For each fixed point $p$ of the action there exist integers $m_1,
\cdots, m_n$ such that the induced linear symplectic $S^1$-action
on the tangent space $T_pM$ is isomorphic to the action on
$(\C^n,\omega_0)$ generated by the moment map
$$
H_0(z_1,\cdots, z_n)=\pi\sum^n_{j=1}m_j|z_j|^2.
$$
The integers $m_1, \cdots, m_n$, uniquely determined up to
permutation, are called the {\it isotropy weights} at $p$.
  An Hamiltonian $S^1$-action on $(M,\omega)$
 is semi-free if and only if the only isotropy weights at every fixed point are
 $\pm 1$.\vspace{2mm}

\noindent{\bf Theorem 6.8}\hspace{2mm}{\it Let $(M,\omega)$ be a
$2n$-dimensional, connected closed symplectic manifold with a
semi-free Hamiltonian circle action with isolated fixed points.
Then
$${\mathcal W}_G(M,\omega)\le C(M,\omega; pt,
PD([\omega]))\le\max H-\min H\leqno(6.14)
$$
for $C=C_{HZ}^{(2)}$, $C_{HZ}^{(2\circ)}$ and any $n\ge 2$, where
$H$ is the associated Hamiltonian function. Moreover, if
 $[\omega]\in H^2(M,\Q)$ and the only isotropy weights at every fixed point is
 $\pm 1$ then
 $${\mathcal W}_G(M,\omega)\ge\frac{\pi}{m}.\leqno(6.15)$$
Here $m>0$ is the smallest integer such that $m[\omega]\in
H^2(M,\Z)$.}\vspace{2mm}

\noindent{\bf Proof.}\hspace{2mm}Following the notations in [Go]
let ${\mathcal S}=\{1,\cdots,n\}$. Each subset $I\subset{\mathcal
S}$ may determine a homology class $A_I\in H_2(M)$ in Proposition
2.11 of [Go] such that $\omega(A_I)=\max H-H(p_{I^c})$ with
$I^c={\mathcal S}\setminus I$. By (14) in Corollary 3.14 of [Go]
one has $x_{\mathcal S}\ast x_I=x_{I^c}\otimes e^{A_I}$. It
follows that Gromov-Witten invariant
 $$ \Psi_{A_I, 0, 3}(pt; PD(x_{\mathcal S}), PD(x_I), PD(x_J))\ne 0\leqno(6.16)
 $$
 for some $J\subset {\mathcal S}$. Note that $x_{\mathcal S}$ is the
  positive generator $H^{2n}(M,\Z)$ (cf. Remark 2.10 in [Go].)
(6.16) shows that $(M,\omega)$ is strong $0$-symplectic uniruled
in the sense of Definition 1.14 in [Lu3, v9]. As in [Lu1] and
[Lu3], using the the reduction formula of the Gromov-Witten
invariants we can also derive from (6.16) that
$$ \Psi_{A_I, 0, 4}(\pi^{-1}(pt); pt, PD([\omega]),\alpha)\ne 0$$
 for some $\alpha\in H_\ast(M,\R)$. So it follows from (12) and Theorem 1.13 in [Lu3, v9] that
\begin{eqnarray*}
{\mathcal W}_G(M,\omega)&&\le C_{HZ}^{(2)}(M,\omega; pt, PD([\omega]))\\
&&\le C_{HZ}^{(2\circ)}(M,\omega; pt, PD([\omega]))\\
&&\le GW_0(M,\omega; pt, PD([\omega]))\\
&&\le\omega(A_I)\\
&&\le\max H-H(p_{I^c})\\
&&\le\max H-\min H.
\end{eqnarray*}
(6.14) is proved.

For the second claim, by Proposition 2.8 in [KaTo] there exists a
symplectic embedding from $(B^{2n}(1),\omega_0)$ to $(M,
m\omega)$. So ${\mathcal W}_G(M, m\omega)\ge\pi$. (6.15) follows.
\hfill$\Box$
\vspace{3mm}


\vspace{2mm}






\begin{thebibliography}{L3}

\bibitem[Ab]{Ab} M. Abreu, {\it K\"ahler geometry of toric manifolds in
                 symplectic coordinates},  Symplectic and contact topology:
                 interactions and perspectives (Toronto, ON/Montreal, QC, 2001),
 Fields Inst. Commun., {\bf 35}, Amer. Math. Soc., Providence,
RI, 2003, pp. 1--24.
\bibitem[Au]{Au} M. Audin,  The topology of torus actions on symplectic manifolds,
                Progress in Mathematics, {\bf 93}: Birkh\"auser, 1991.
\bibitem[Ba1]{Ba1} V. V. Batyrev, {\it Quantum cohomology rings of toric manifolds},
                Ast\'erisque, {\bf 218}(1993), 9-34.
\bibitem[Ba2]{Ba2} V. V. Batyrev, {\it On the classification of smooth projective
                  toric varieties},  J. Algebraic Geometry, {\bf 3}(1994), 493-535.
\bibitem[Ba3]{Ba3} V. V. Batyrev, {\it On the classification of toric Fano $4$-folds},
                   J.Math.Sciences, {\bf 94}(1999), 1021-1050.
\bibitem[Bi]{Bi} P. Biran, {\it From symplectic packing to algebraic geometry and back},
                   European Congress of Mathematics, {\bf Vol}. II(Barcelona,2000),
                   Prog. Math.,202: Birkh\"auser, 2001, pp. 507-524.
\bibitem[BiCi]{BiCi} P. Biran and K. Cieliebak, {\it Symplectic topology on
                     subcritical manifolds},  Comm. Math. Helv.,
                     {\bf 76}(2001), no.4, 712-753.
\bibitem[CdFKM]{CdFKM} P. Candelas, X. de la Ossa, A. Font, S. Katz, D. Morrison,
                   {\it Mirror symmetry for two parameter models I},
                     Mirror symmetry, II,  AMS/IP Stud. Adv. Math., {\bf 1},
                     Amer. Math. Soc., Providence, RI, 1997, pp.
                     483--543.
\bibitem[CiS]{CiS}  K. Cieliebak and D.A. Salamon, {\it Wall crossing for symplectic
                     vortices and quantum cohomology},  math.SG/0209170.
\bibitem[Del]{Del} T. Delzant, {\it Hamiltoniens p\'eriodiques et image convexe de
                  l'application moment},  Bull. Soc. Math. France,
                 {\bf 116}(1988), 315-339.
\bibitem[Dem]{Dem} J.-P. Demailly, {\it $L^2$-vanishing theorems for positive line
                    bundles and adjunction theory}, In  Transcendental
                   methods in Algebraic Geometry (ed. F.Catanese and C.Ciliberto),
                   Lect. Notes Math. {\bf 1646}, Springer-Verlag, 1992, pp. 1-97.
\bibitem[Ew]{Ew} G. Ewald,  Combinatorial combinatorial convexity and algebraic
                   geometry,  Graduate Texts in Mathematics {\bf 168}, Springer, 1996.
\bibitem[Fu]{Fu} W. Fulton,  Introduction to Toric Varieties,
                 Annals of Mathematics Studies {\bf 131}, Princeton University Press, 1993.
\bibitem[Ga]{Ga} A. Gathmann, {\it Gromov-Witten invariants of blow-ups},
                   J. Algebraic Geom., {\bf 10}(2001), no. 3, 399--432.
\bibitem[Gin]{Gin} V. Ginzburg, {\it The Weinstein conjecture and the theorems of nearby
                    and almost existence},  The breadth of symplectic and Poisson geometry,
                     139--172, Progr. Math., {\bf 232},  Birkh\"user Boston, Boston, MA, 2005.
\bibitem[Giv]{Giv} A. Givental, {\it A mirror theorem for toric complete intersections},
                      Topological field theory, primitive forms and related topics (Kyoto, 1996),
                      141--175, Progr. Math., {\bf 160},  Birkh\"user Boston, Boston, MA,
                      1998.
\bibitem[Go]{Go} E. Gonzalez, {\it Quantum cohomology and $S^1$-action with isolated fixed points},
                  math.SG/0310114.
\bibitem[Gr]{Gr} M.~Gromov, {\it Pseudoholomorphic curves in symplectic manifolds},
                 Invent. Math., {\bf 82} (1985), 307-347.
\bibitem[Gu1]{Gu1} V. Guillemin, {\it K\"ahler structures on toric varieties},
                   J. Diff. Geom., {\bf 40} (1994), 285-309.
\bibitem[Gu2]{Gu2} V. Guillemin,  Moment maps and combinatorial invariants of
                    Hamiltonian $\T^n$-spaces,  Progress in Mathematics,
                  {\bf 122}: Birkh\"auser, 1994.
\bibitem[HaKn]{HaKn} J.-C. Hausmann and A. Knutson, {\it The cohomology ring of polygon spaces},
                    Ann. Inst. Fourier, Grenoble, {\bf
                   48}(1998), 281-321.
\bibitem[Hu]{Hu} Jianxun Hu, {\it Gromov-Witten invariants of blow-ups along points and
                  curves},  Math. Z., {\bf 233}(2000), no. 4, 709--739.
\bibitem[HZ]{HZ} H.~Hofer and E.~Zehnder,  Symplectic Invariants and
                  Hamiltonian Dynamics, Birkhauser, Boston, MA (1994).
\bibitem[Ka]{Ka} Y. Karshon, {\it Appendix to} [McPo],  Inven. Math., {\bf 115}(1994),
                  431-434.
\bibitem[KaTo]{KaTo} Y. Karshon and S. Tolman, {\it The Gromov width of
                    complex Grassmannians},  math.SG/0405391.
\bibitem[Ko]{Ko} J. Koll\'ar, {\it Low Degree Polynomial Equations: Arithmetric, Geometry and
                 Topology},   Progress in Mathematics,
                  122: Birkh\"auser, 1994, pp.255-288.
\bibitem[KoMor]{KoMor} J. Koll\'ar and S. Mori,  Birational Geometry of
Algebraic Varieties, Cambridge tracts in Mathematics, 134:
Cambridge University Press, 1998.

\bibitem[Kr]{Kr} A. Kresch, {\it Gromov-Witten invariants of a class of toric varietes},
                    Michigan Math. J., {\bf 48} (2000), 369--391.
\bibitem[Lu1]{Lu1} G. C. Lu, {\it The Weinstein conjecture in the
                       uniruled manifolds},  Math. Res. Lett., {\bf 7}(2000),383-387.
\bibitem[Lu2]{Lu2} G. C. Lu, {\it Symplectic capacities of toric manifolds and
                   combinatorial inequalities},  C. R. Acad. Sci. Paris, Ser. I,
                   {\bf 334}(2002), 889-892.
\bibitem[Lu3]{Lu3} G. C. Lu,  {\it Gromov-Witten invariants and pseudo symplectic capacities},
                    math.SG/0103195, v6, 6 September 2001, and v9, 3 December 2004.
                    to appear in  Israel Journal of Mathematics.
\bibitem[Mc]{Mc}  D. McDuff, {\it Quantum homology of fibrations over $S^2$},
                   International Journal of mathematics, {\bf
                   11}(2000), 665-721.
\bibitem[McPo]{McPo} D. McDuff and L. Polterovich, {\it Symplectic packings and
                      algebraic geometry},  Invent. Math., {\bf 115}(1994),
                      405-425.
\bibitem[Mor1]{Mor1} S. Mori, {\it Projective manifolds with ample tangent bundles},
                      Ann. Math., {\bf 110}(1975), 593-606.
\bibitem[Mor2]{Mor2} S. Mori, {\it An email communication}.
\bibitem[Oda]{Oda} T. Oda,  Convex Bodies and Algebraic Geometry,
                  Springer-Verlag, 1988.
\bibitem[Sa]{Sa} H. Sato, {\it Toward the classification of higher-dimensional toric Fano varieties},
                    Tohoku Math. J., {\bf 52}(2000), no. 3, 383--413.
\bibitem[Sch]{Sch} F. Schlenk, {\it On symplectic folding},  preprint,
                   math.SG/9903086,  March 1999.
\bibitem[Sik]{Sik} J.C. Sikorav, {\it Rigidit\'e symplectique dans le cotangent de
                  $\T^n$},  Duke Mathematical Journal, {\bf 59}(1989), 227-231.
\bibitem[Sp]{Sp} H.Spielberg, {\it The Gromov-Witten invariants of  symplectic toric
                  manifolds, and their quantum cohomology ring},
                   C. R. Acad. Sci. Paris, Ser. I, {\bf 329}(1999), 699-704.
\bibitem[Tr]{Tr} L. Traynor, {\it Symplectic packing constructions},  J. Diff.
                  Geom., {\bf 41}(1995), 735-751.
\bibitem[Wi]{Wi} J.A.Wi\'sniewski, {\it Toric Mori theory and Fano manifolds},
                 S\'eminaires \& Congr\`es, {\bf 6}(2002), 249-272.
\end{thebibliography}
\end{document}